\documentclass[a4paper]{article}
\usepackage{amssymb}
\usepackage{amsmath}
\usepackage{mathrsfs}
\usepackage{amsfonts}
\usepackage{color}
\usepackage{vmargin}
\usepackage{amsthm}
\usepackage{comment}

\long\def\symbolfootnote[#1]#2{\begingroup%
\def\thefootnote{\fnsymbol{footnote}}\footnote[#1]{#2}\endgroup}

\setmarginsrb{20mm}{20mm}{20mm}{20mm}{10mm}{10mm}{10mm}{10mm}

\numberwithin{komcounter}{section}

\newtheorem{theorem}{Theorem}

\newtheorem{proposition}{Proposition}[section]

\newtheorem{corollary}[proposition]{Corollary}
\newtheorem{definition}[proposition]{Definition}
\newtheorem{example}[proposition]{Example}

\newtheorem{lemma}[proposition]{Lemma}

\newtheorem{remark}[proposition]{Remark}

\title{Radial reductions of variable exponent Laplacian equations on warped products with compact symmetry}
\author{Jos\'e Luis D\'iaz Palencia\footnote{joseluis.diaz.p@udima.es }, ~ Daniel J. Pons \footnote{dpons@unab.cl ; pons.dan@gmail.com } ~ and ~ Enrique G. Reyes
\footnote{enrique.reyes@usach.cl ; e\_g\_reyes@yahoo.ca} \\
\\}

\begin{document}

\maketitle

\begin{abstract}
We study nonlinear elliptic equations of variable exponent type,
\[
-\Delta_{q(x)}u+|u|^{q(x)-2}u=f(x,u),
\]
on non-compact Riemannian manifolds with symmetry. We consider complete warped products
\[
M=\mathcal R\times_\phi Y,\qquad
g_\phi=dr^2\oplus\phi(r)^2g_Y,
\]
where $\mathcal R$ is a non-compact interval, possibly with boundary, and $Y$ is a compact homogeneous manifold. 
We establish sufficient conditions on the warping function ensuring a Ricci lower bound, a uniform non-collapsing volume 
condition, and the orbit-growth hypothesis of the abstract existence theorem proved in \cite{JFA}. 
We also obtain comparison bounds for the admissible warping functions.

We analyze the corresponding radial equations for the admissible exponential and linear models $\phi(r)=e^{ar}$ and $\phi(r)=ar$, with $r\geq r_0>0$. 
For a constant diffusion exponent $q=k$, we consider explicit decaying profiles $u=\phi^{-\sigma}$, which belong to $L_1^k(M)$ provided $\sigma k>\dim Y$ 
in the exponential case and $\sigma k>\dim Y+1$ in the linear case. Under additional explicit parameter restrictions, these profiles solve a resonant equation 
with a positive bounded coefficient, a construction distinct from the superlinear regime of the abstract existence theorem. Finally, we establish logarithmic 
reductions for constant and variable radial exponents, and justify the passage from radial identities to weak formulations against arbitrary compactly supported interior test functions.
\end{abstract}

\paragraph{Keywords:} Variable exponent Laplacian; warped product manifolds; radial reductions; compact group actions; Ricci curvature bounds; weak solutions; explicit radial solutions.

\smallskip

\paragraph{2020 Mathematics Subject Classification:} Primary: 35J92. Secondary: 53C21, 58J05, 46E35, 35A24.

\section{Introduction}

\begingroup

The purpose of this paper is to study nonlinear elliptic equations of
variable exponent type on non-compact Riemannian manifolds with symmetry.
The model equation is
\begin{equation}\label{intro:eq}
-\Delta_{q(x)}u+|u|^{q(x)-2}u=f(x,u)\; ,
\end{equation}
where
\[
\Delta_{q(x)}u
=\operatorname{div}_g\left(|\nabla_g u|^{q(x)-2}\nabla_g u\right)
\]
is the variable exponent Laplacian associated with the Riemannian metric $g$.
Allowing the exponent to depend on position describes a spatially varying growth
law and removes the fixed homogeneity of the classical $k$-Laplacian.  
Such operators arise naturally in models of electrorheological fluids,
where the constitutive response depends on an applied electric field;
see R\r{u}\v{z}i\v{c}ka~\cite{Ruzicka}.
Another relevant application arises in electrothermal systems.
Bul\'{\i}\v{c}ek, Glitzky and Liero~\cite{BulicekGlitzkyLiero}
established the existence of weak solutions for coupled elliptic
systems with $p(x)$-Laplacian-like constitutive laws, allowing
discontinuous spatial exponents. Their results demonstrate the
relevance of variable-growth operators in heterogeneous materials,
where local transport properties may vary abruptly across
different regions. 

The functional framework is provided by variable exponent Lebesgue and Sobolev
spaces. Foundational results on their norm, duality and embedding properties were
developed by Kov\'a\v{c}ik and R\'akosn\'ik~\cite{KovacikRakosnik} and by Fan and
Zhao~\cite{FanZhao}; a systematic account is given in the monograph of Diening,
Harjulehto, H\"ast\"o and R\r{u}\v{z}i\v{c}ka~\cite{Diening}.
For nonlinear equations, the variational treatment of the $p(x)$-Laplacian
Dirichlet problem by Fan and Zhang~\cite{FanZhang} illustrates how this framework
extends familiar existence methods beyond constant exponents.

Passing from Euclidean domains to complete non-compact manifolds introduces
geometric restrictions into the Sobolev embedding theory. Curvature bounds and control
of local volume prevent degeneracies that cannot be detected from a single
coordinate chart. The classical Sobolev theory on Riemannian manifolds is
developed in Hebey~\cite{Hebey}, while the relevant comparison geometry is
discussed in Petersen~\cite{Petersen}. The dependence of analysis on both the
metric and the measure is also central to the geometric viewpoint presented in
Villani~\cite{Villani}. In our setting, a lower Ricci bound and a uniform
positive lower bound for volumes of unit balls provide the geometric background
for the existence theorem recalled below. These assumptions control local
geometry uniformly, but additional information is needed to recover global
compactness at infinity, see Theorem~\ref{reac-diff}.

The obstruction is familiar: bounded sequences may carry their mass farther and
farther away while retaining bounded energy. Concentration-compactness provides
a general framework for this phenomenon~\cite{LionsCC}. Symmetry can restrict
such escape, as illustrated by the radial compactness underlying Strauss's
work~\cite{Strauss} and the more general symmetry results of
Lions~\cite{LionsSymmetry}. On manifolds, the influence of isometry group actions on
Sobolev embeddings was studied by Hebey and Vaugon~\cite{HebeyVaugon}.
For non-compact homogeneous manifolds, Skrzypczak and
Tintarev~\cite{SkrzypczakTintarev} relate compactness of invariant subspaces to
the geometry of the group orbits near infinity. These results motivate a
concrete question: which warped product geometries make the compactness
mechanism effective for equations with non-standard diffusive growth?

Our principal analytical input is a theorem in ~\cite{JFA}, recalled as Theorem~\ref{reac-diff}.
It gives a non-trivial invariant weak solution of~\eqref{intro:eq} under
subcriticality and regularity assumptions on the exponents, suitable growth
conditions on $f$, and an expansion condition on the orbits of a compact group
of isometries. The superlinearity assumption belongs to the variational
tradition of Ambrosetti and Rabinowitz~\cite{AmbrosettiRabinowitz}.
The orbit condition is expressed through the number $M(x,R)$ of disjoint
ambient balls of a fixed radius whose centres lie in the orbit of $x$.
Requiring this number to diverge at infinity has a direct interpretation:
an invariant concentration situated far away would have to be replicated
across an increasing number of disjoint regions, which is incompatible with
a fixed bound on its total mass. The purpose here is to examine the geometric
content of this hypothesis and the radial equations associated with it.

We focus on warped products
\[
(M,g_\phi)=
(\mathcal R\times Y,\,dr^2\oplus\phi(r)^2g_Y),
\]
where $\mathcal R$ is a complete non-compact interval, possibly with boundary,
and $(Y,g_Y)$ is a compact homogeneous manifold for a compact Lie group $H$.
The action is $h\cdot(r,y)=(r,h\cdot y)$, so that the orbits are precisely the
fibres $\{r\}\times Y$. The warped product construction originates in the
classical work of Bishop and O'Neill~\cite{BishopONeill}; its curvature and
submersion geometry are treated in Besse~\cite{Besse} and Gromoll and
Walschap~\cite{Gromoll}. If $m= n-1=\dim Y$, then
\[
dV_{g_\phi}=\phi(r)^m\,dr\wedge dV_{g_Y},
\qquad
\Delta_{\phi,q}u
=\frac{1}{\phi(r)^m}
\left(\phi(r)^m|u_r|^{q(r)-2}u_r\right)_r
\]
for radial functions. Consequently, the same warping function determines the
curvature, the growth of the symmetry orbits, the radial drift, and the weight
against which a candidate solution must be integrable. This makes warped
products a useful setting in which to compare geometric hypotheses with
explicit analytical calculations.

The first direction of the paper is geometric. For compact homogeneous Einstein
fibres, Theorem~\ref{Warp-1D} gives sufficient one-dimensional inequalities for
a Ricci lower bound. Corollary~\ref{coro-warp} describes the associated
comparison functions, and Section~\ref{models} examines the volume and orbital
conditions required by Theorem~\ref{reac-diff}.
The choice of a canonical model must respect compact homogeneity of $Y$.
In the Einstein classification discussed in~\cite{Besse,Solitons}, the
$\cosh$ branch with negative Einstein constant on the fibre is incompatible
with that requirement; see Remarks~\ref{H-spaces} and~\ref{not-cosh}.
This exclusion concerns that Einstein branch. The explicit constructions in
Section~\ref{sec:radial} use the exponential and linear warpings
$\phi(r)=e^{ar}$ and $\phi(r)=ar$, $a>0$, on expanding ends
$[r_0,\infty)$, with $r_0>0$. Keeping the interval explicit is essential:
the geometry and integrability at one end do not determine their behaviour on
a different global extension.

The second direction concerns the reduced equation and the meaning of its
solutions. Transitivity of $H$ on $Y$ makes invariant functions radial.
For constant $q=k>2$, a logarithmic substitution yields equivalent flux and
divergence identities, as well as explicit constant and $k$-harmonic profiles.
The corresponding substitution for $q=q(r)$ exposes the additional dependence
on $q_r$ and provides a structural identity for the variable exponent problem.
The relation between invariant variations and unrestricted weak formulations
is consistent with Palais's principle of symmetric
criticality~\cite{Palais}; here the relevant implication follows from
\cite[Lemma 6.1]{JFA} and is recorded in
Corollary~\ref{cor:unrestricted-weak-formulation}.
We also justify directly the passage from a radial distributional identity to
a weak identity on the full warped product by integration over the compact
fibre.

Particular attention is paid to decay in the weighted Sobolev space.
Non-zero constant profiles fail to be integrable on the infinite-volume models,
whereas Proposition~\ref{prop:decaying-manufactured-profile} constructs
\[
u(r)=\phi(r)^{-\sigma}
=\begin{cases}
e^{-\sigma ar},&\phi(r)=e^{ar},\\
(ar)^{-\sigma},&\phi(r)=ar.
\end{cases}
\]
On $[r_0,\infty)\times_\phi Y$, the respective conditions (recall $m = n-1$),
$\sigma k>m$ and $\sigma k>m+1$ ensure that both $u$ and its gradient are
$k$-integrable. With the parameter restrictions stated in that proposition,
the profiles solve
\[
-\Delta_{\phi,k}u+u^{k-1}=\gamma(r)u^{k-1}
\]
for an explicit positive bounded coefficient $\gamma$.
This is the resonant case $p=k$, and is therefore distinct from the
superlinear regime $p>q=k$ required by Theorem~\ref{reac-diff}.
The profiles exhibit how the geometry determines admissible decay and how an
explicit strong solution can belong to the variational space; they are not
claimed to realize every hypothesis of the abstract theorem.
On a half-line model, their weak formulation uses test functions compactly
supported in the interior and does not impose a boundary condition at $r_0$.

The paper is organized as follows. Section~\ref{sec:prelim} collects the warped
product curvature formulas, the variable exponent spaces, and the abstract
existence theorem. Section~\ref{sec:geometry} studies invariant functions,
the unrestricted weak formulation, and geometric criteria for warped products
with one-dimensional base. Section~\ref{sec:radial} develops the constant and
variable exponent logarithmic reductions, discusses the integrability of the
explicit radial families, constructs decaying profiles on the exponential and
linear models, and establishes their weak formulation on the full product.

\par
\endgroup

\section{Preliminaries}\label{sec:prelim}

Given  a Riemannian $n$-dimensional manifold $(M,g)$, we denote by $TM$ its tangent bundle, and by  $T^* M$ its cotangent bundle, so that $g$ is a positive definite symmetric section of 
$T^*M \otimes T^*M$.  Two classical differential operators associated to $g$ are $\nabla$ and $\text{div}$, the  
Levi-Civita connection and divergence, respectively. We recall that $\nabla$ is defined uniquely via the Koszul
formula
\[
2g(\nabla_\tau\varsigma,\eta) = 
\tau\bigl(g(\varsigma,\eta)\bigr) + \varsigma\bigl(g(\tau,\eta)\bigr) - \eta\bigl(g(\tau,\varsigma)\bigr) - g(\tau,[\varsigma,\eta]) + g(\varsigma,[\eta,\tau])  + g(\eta,[\tau,\varsigma])\; ,
\]
in which $\tau$, $\varsigma$ and $\eta$ are sections of $TM$, and that (assuming that $M$ is orientable and equipped with the standard Riemannian measure, as we will do from now on) 
$\text{div}$ is determined by the equation
$$
L_\tau dV_g = \text{div}(\tau) dV_g \; ,
$$
in which $L$ denotes Lie derivative, $\tau$ is a section of $TM$ and $dV_g$ is the Riemannian volume $n$-form of $(M,g)$. Let $u$ be a $C^2$ 
function on $M$; then:
\begin{itemize}
\item Its Hessian, regarded as a symmetric section of $T^*M \otimes T^*M$, is given by $\text{Hess}\ u := \nabla \nabla u$, so that 
\[
\operatorname{Hess} u(\tau,\varsigma)
=
\tau(\varsigma u)-(\nabla_\tau\varsigma)u.
\]
Equivalently,
\[
\operatorname{Hess} u(\tau,\varsigma)
=
g(\nabla_\tau\operatorname{grad} u,\varsigma).
\]
whenever $\tau$ and $\varsigma$ are sections of $TM$.
\item Its Laplacian, regarded as a continuous function on $M$, is given by
\[
\Delta u
=
\operatorname{tr}(\operatorname{Hess} u)
=
\operatorname{div}(\operatorname{grad} u).
\]
\end{itemize} 
Sometimes we use the subscript $g$ if the metric used  is not obvious  from the context.


\subsection{Warped products}

We consider an $n$-dimensional connected and orientable manifold $(M,g_{\phi})$ of the form
$M = X \times Y$,   in which  $(Y,g_Y)$ is a $m$ dimensional manifold, and $(X,g_X)$ is an $(n-m)$-dimensional 
manifold, endowed with the warped product metric  
\begin{equation}     \label{warped}
g_{\phi}(x,y) = g_X(x) \oplus \phi(x)^2 g_Y(y) 
\end{equation}
for $(x,y) \in M$, and where $\phi : X \to \mathbb{R}_+$ is a smooth given function.

The  canonical volume $n$-form on $(M,g_{\phi})$  is decomposed as
$$d V_{g_{\phi}}(x,y) = d V_{g_X}(x) \wedge \phi(x)^{m}  dV_{g_{Y}}(y) = 
\phi(x)^{m} d V_{g_X}(x) \wedge d V_{g_Y}(y)\ ,$$ where
$dV_{g_{X}}$ is the  $(n-m)$-volume form on
$(X, g_X)$, and  $dV_{g_{Y}}$ is the  $m$-volume form on
$(Y, g_Y)$. 

Every vector field in $M$, $\tau$, say, is decomposed into its \textit{horizontal} part $\tau_X$  in $T X$, and 
its \textit{vertical} part  $\tau_Y$ in $T Y$, so that
$\tau = \tau_X \oplus \tau_Y $, i.e. $ TM = \pi_X^* T X \oplus \pi_Y^* TY$, in which $\pi_X$ and $\pi_Y$ are standard projections. We will write simply $TM = TX \oplus TY$. 
This decomposition extends to the tensor products of $T M$ and $T^{*}M$. 

In particular, we can make the identification
\begin{equation}   \label{tensor}
 T^{*}M \otimes T^{*}M  = ( T^{*}X \otimes T^{*}X ) \oplus  ( T^{*}X \otimes T^{*}Y ) \oplus  
( T^{*}Y \otimes T^{*}X ) \oplus  ( T^{*}Y \otimes T^{*}Y ) \ , 
\end{equation}
so that for warped product metrics of the form (\ref{warped}), we obtain
$$  g_{\phi}(\tau, \varsigma)  = g_{X}(\tau_X , \varsigma_X) + \phi^{2} g_{Y}(\tau_Y , \varsigma_Y)  \ .$$
The decomposition for the Ricci tensor of  $g_{\phi}$, using (\ref{tensor}), gives
\begin{equation}
\text{Rc}(g_{\phi})(\tau,\varsigma) = \text{Rc}(g_{\phi})(\tau_X,\varsigma_X)  +   \text{Rc}(g_{\phi})(\tau_X,
\varsigma_Y)  + \text{Rc}(g_{\phi})(\tau_Y,\varsigma_X) + \text{Rc}(g_{\phi})(\tau_Y,\varsigma_Y) \ .
\end{equation}

\begin{lemma}  \label{Ricci-warp}
Using the previous notation:
\begin{enumerate}
\item  \label{hor} $  \text{Rc}(g_{\phi})(\tau_X,\varsigma_X) =  \text{Rc}(g_X)(\tau_X,\varsigma_X) - 
\frac{m}{\phi}\  \text{Hess}\ \phi  (\tau_X,\varsigma_X) $
\item $  \text{Rc}(g_{\phi})(\tau_X,\varsigma_Y)  = 0 =   \text{Rc}(g_{\phi})(\tau_Y,\varsigma_X)  \  , $
and 
\item \label{vert} $\displaystyle  \text{Rc}(g_{\phi})(\tau_Y,\varsigma_Y) =  \text{Rc}(g_Y)(\tau_Y,\varsigma_Y) -
(\ \frac{1}{m}\ \frac{\Delta_{g_X} \phi^m}{\phi^m} \ )   \phi^2 g_Y(\tau_Y,\varsigma_Y) \  .  $
 \end{enumerate}
\end{lemma} 

\begin{proof}
This proof follows from classical identities, see Proposition 9.106  in \cite{Besse}, or Corollary 2.2.2   
in \cite{Gromoll}, in addition to
$$  m\ \frac{\text{Hess}\ \phi}{\phi}  = 
\text{Hess}\ (m \log \phi) + \frac{1}{m} \nabla (m \log \phi) \otimes \nabla (m \log \phi)   $$
in $T^{*}X \otimes T^{*}X $, and 
$$ \frac{1}{m}\ (\Delta_{g_X} \phi^m ) / \phi^m  =   
(\Delta_{g_X} \phi) / \phi + (m-1) \vert \nabla \log \phi \vert^2_{g_X} \ ,  $$
if $\phi$ is of type $C^2$. 
\end{proof}


\subsection{Notation and framework from \cite{JFA}}

We recall some notation and framework from \cite{JFA}. Hereafter, $B_r(p)$ denotes the open ball of radius $r$ centered at $p$, $p \in M$.

\begin{definition}    \label{B-vol}
We say that the Riemannian $n$-manifold $(M,g)$ has the property $\mathcal{B}_{\text{vol}}(\lambda,\vartheta)$ if:
\begin{itemize}
\item $Rc(g) \geq \lambda (n-1) g$ for some $\lambda \in \mathbb{R}$, and
\item There exists some $\vartheta > 0$ such that $| B_1(x)  |_g \geq \vartheta$ for every $x$ in $M$. 
\end{itemize}
\end{definition}


\begin{definition}
We denote by $\text{Iso}(M,g)$ the finite-dimensional Lie group of isometries of $(M,g)$. 
Let $H$ be a Lie subgroup of $\text{Iso}(M,g)$, and fix $x$ in $M$ and $R > 0$. 
The number $M(x,R) \in \mathbb{N} \cup \{+\infty \} $ is the lowest upper bound for the number of disjoint $R$-balls with centers in the $H$-orbit of $x$, that is,
\[
M(x,R):=
\sup\left\{
N\in\mathbb N:\ \exists x_1,\ldots,x_N\in Hx
\text{ such that }
B_R(x_i)\cap B_R(x_j)=\varnothing
\text{ for }i\neq j
\right\}.
\]\end{definition}

In this definition we use the standard fact, due to Myers--Steenrod, that the isometry
group of a connected smooth Riemannian manifold is a Lie group.
The abstract result that we use in this paper is Theorem \ref{reac-diff} below. We need to introduce the following class of functional spaces.

\begin{definition}
Let $(M,g)$ be a smooth Riemannian manifold with Riemannian measure
$d V_g$, and let
\[
q:M\to (1,\infty)
\]
be a measurable function. We set
\[
q^-:=\operatorname*{ess\,inf}_{x\in M}q(x),
\qquad
q^+:=\operatorname*{ess\,sup}_{x\in M}q(x),
\]
and assume throughout that
\[
1<q^-\leq q^+<\infty .
\]
The variable exponent Lebesgue space $L^{q(\cdot)}(M)$ is defined by
\[
L^{q(\cdot)}(M)
:=
\left\{
u:M\to\mathbb R \ \text{measurable} \;:\;
\int_M |u(x)|^{q(x)}\,d V_g(x)<\infty
\right\}.
\]
It is equipped with the norm
\[
\|u\|_{L^{q(\cdot)}(M)}
:=
\inf\left\{
\lambda>0 \;:\;
\int_M
\left|\frac{u(x)}{\lambda}\right|^{q(x)}
\,d V_g(x)
\leq 1
\right\}.
\]
The corresponding variable exponent Sobolev space is
\[
L_1^{q(\cdot)}(M)
:=
\left\{
u\in L^{q(\cdot)}(M)
\;:\;
|\nabla_g u|\in L^{q(\cdot)}(M)
\right\},
\]
where $\nabla_g u$ denotes the weak Riemannian gradient of $u$. We endow
$L_1^{q(\cdot)}(M)$ with the norm
\[
\|u\|_{L_1^{q(\cdot)}(M)}
:=
\|u\|_{L^{q(\cdot)}(M)}
+
\|\nabla_g u\|_{L^{q(\cdot)}(M)} .
\]
\end{definition}

\begin{definition} \label{H-inv}
Whenever $ X $ is a subset of the real valued functions defined on $ M $, and $ H $ is a group acting on $ M $, 
we denote by $ X_H $ the subset of $ M $ made up of the $ H $-invariant functions in $ X $, so that
\[   X_H := X \cap \{\  u : \ h^{\ast} u = u \ \text{whenever}\ h \in H\ \}.    \]
\end{definition}

\begin{definition}
We use $\mathcal{D}(M)$ to denote the space of test functions of compact support on $M$. We consider the equation
\begin{equation} \label{eqExample}
-\Delta_{q(x)} u(x) + |u(x)|^{q(x)-2}u(x) = f(x,u(x)) \ ,
\end{equation}
in which
\begin{equation*}
\Delta_{q(x)} u := div_g  \left( |\nabla_g u|^{q(x)-2} \nabla_g u  \right) 
\end{equation*}
is the formal $q(x)$-Laplacian on $(M,g)$. \\ 

\noindent We say that $u\in L^{q(\cdot)}_1 (M)$ is a weak solution of $(\ref{eqExample})$ if for every $\psi \in \mathcal{D}(M)$ we have
\begin{equation}   \label{weak-sol}
 \int_M \left(|\nabla_g u(x)|^{q(x)-2} g(\nabla_g u(x) , \nabla_g \psi (x))  + |u(x)|^{q(x)-2}u(x)\psi (x) - f(x,u(x))\psi(x)\right) dV_g (x) =0\; ,
\end{equation}
and we shall call $u\in L_1^{q(\cdot)}(M)$ an $H$-weak solution if the weak
identity holds for every $\psi\in\mathcal D_H(M)$, in which $\mathcal{D}_H(M)$-see Definition \ref{H-inv} above-is the space of $H$-invariant functions in $\mathcal{D}(M)$. 
\end{definition}


\smallskip

We introduce the following notation:  $\mathcal{P}(M)$ is the set of variable Sobolev exponents, that is, the set of bounded measurable functions from 
$M$ to $[1,\infty)$; $q^*(x) = n q(x)/(n-q(x))$ whenever $q^+<n$,  and for variable exponents $r,s\in\mathcal P(M)$, we write $r\ll s$ if $\operatorname*{ess\,inf}_{x\in M}(s(x)-r(x))>0$. 
Thus $q\ll p\ll q^*$ means that
\[
\operatorname*{ess\,inf}_{x\in M}(p(x)-q(x))>0,
\qquad
\operatorname*{ess\,inf}_{x\in M}(q^*(x)-p(x))>0. \]

In \cite{JFA} the following theorem was proved:


\begin{theorem} (Theorem 6.1 in \cite{JFA})  \label{reac-diff}

\textbf{Hypotheses}: 
\begin{enumerate}
\item $H$  is a non-trivial compact Lie subgroup of $\text{Iso}(M,g)$. 
\item \label{bounded-geo} The $n$-manifold $(M,g)$ is complete and it has property $\mathcal{B}_{\text{vol}}(\lambda,\vartheta)$.
\item \label{M-balls} If $R>0$ and $o$ is any fixed point in $M$, then
$$ \displaystyle \lim_{ T \rightarrow \infty} \, \inf_{x \in M \setminus B_T(o)} M(x,R) = \infty .$$
\item We take $q$ and $p$ in $\mathcal{P}(M)$; we assume that they are $H$-invariant, that $q$ is uniformly continuous, and that 
\begin{equation*}
1 < q^-\leq q^+ < n\ ,\  q \ll p \ll q^*\ ,\  \text{and}\  q^+ < p^-\ .
\end{equation*}
\item  \label{f-term}  $f: M \times \mathbb{R} \rightarrow \mathbb{R}$ is a continuous function satisfying:   
\begin{itemize}
\item [($i$)] For each $t \in \mathbb{R}$ the function $f( \cdot ,t) : M \rightarrow \mathbb{R}$ is $H$-invariant\ ,
\item[($ii$)] There exists $c_1 > 0$ such that for each $t \in \mathbb{R}$ and all $x \in M$ we have the bound
$$ |f(x,t)| \leq c_1  |t|^{p(x) -1}\; , $$ 
\item[($iii$)] There exists $\theta > q^+$ and $A > 0$ such that for each $|t| > A $ and $x \in M$ we have  
\begin{equation*}
0 < F(x,t) \leq f(x,t)\ t / \theta\ ,
\end{equation*} 
where $F(x,t) := \int_0^t f(x,s) \, ds$.
\end{itemize} 
\end{enumerate}

\textbf{Conclusion:} Equation $(\ref{eqExample})$ has an $H$-invariant  non-trivial $H$-weak solution in 
$L_1^{q(\cdot)}(M)$, that is there exists some $ u \in L_{1,H}^{q(\cdot)}(M) $ such that $(\ref{weak-sol})$ holds for every  $\psi$ in $\mathcal{D}_H(M)$.
\end{theorem}


\section{$H$-invariant solutions on warped products with one-dimensional base}\label{sec:geometry}

\subsection{$H$-invariant functions}
As noticed in \cite[Remark 3]{JFA}, applications of  Theorem \ref{reac-diff} arise naturally in the following picture: Consider $n$ dimensional manifolds $(M,g_{\phi})$ of the form
$M = X \times Y$, with warped product metrics as in (\ref{warped}), namely  
$$ g_{\phi}(x,y) = g_X(x) \oplus \phi(x)^2 g_Y(y)  , $$
where $(X,g_X)$  is a complete non-compact $(n-k)$-manifold,  and $(Y,g_Y)$ is a  $k$-manifold that is  a homogeneous space for a compact Lie group $H$ that acts isometrically on $(M,g_{\phi})$ 
and leaves the elements of $X$ fixed (in particular, $Y$ is a compact manifold). Explicitly, whenever $h \in H$ and $(x,y) \in M$, we have that $h(x,y) = (x, h y)$, so that the orbit $H(x,y)$ is 
$\{  x \} \times Y$, with 
$$h^{\ast} g_{\phi} = g_X \oplus \phi(x)^2 h^{\ast} g_Y $$
and $h^{\ast} g_Y = g_Y$.  In other words, $\pi : (M,g_{\phi}) \to (X, g_X)$ is a Riemannian submersion, and the 
fiber $\pi^{-1}(x)$ over $x \in X$ is isometric to $(Y, \phi(x)^2 g_Y)$,  with
$H$ acting transitively and isometrically in each fiber, see \cite{Gromoll}. \\

\begin{remark}     \label{H-spaces}
A discussion on homogeneous spaces, that we follow, is given on Chapter 7 of \cite{Besse}. Consider
$(Y, g_Y)$, a homogeneous space for $H$.  An important subgroup of $H$ is the isotropy subgroup, say $K$, so that $H$ acts effectively on $H/K$, and we have the diffeomorphism $Y \equiv H/K$. 

In this identification, the Lie algebra of $H$, say $\mathfrak{H}$, has a direct sum decomposition $\mathfrak{H} = \mathfrak{K} \oplus \mathfrak{P}$: Given $y$ in $Y$,  $\mathfrak{P}$ can be identified with  $T_yY$, the non-trivial Killing vector fields for $(Y,g_Y)$ at $y$, and $\mathfrak{K}$ is the Lie algebra of $K$. 
 
Let $\text{Ad} : H \to \text{GL}(\mathfrak{H})$ denote the adjoint representation of $H$; using the previous identifications, we have a natural representation $\text{Ad} : K \to  \text{GL}(\mathfrak{P})$, called the linear isotropy representation  of $K$ in $\text{GL}(\mathfrak{P})$. 
 
 A result attributed to J. Wolf, says that if the linear isotropy representation of $K$ in $\text{GL}(\mathfrak{P})$
 is irreducible, then $Y$ admits an $H$-invariant Einstein metric, that is unique up to dilations. For instance, all irreducible symmetric spaces satisfy this condition. 
 
Finally, we mention that if $(Y,g_Y)$ is an Einstein homogeneous space for a compact Lie group $H$ then the Einstein constant is non-negative. 
 \end{remark}

Theorem~\ref{reac-diff} yields a critical point in the $H$-invariant
Sobolev space, and the weak formulation is initially obtained only
against $H$-invariant test functions, that is  $(\ref{weak-sol})$ holds for every $\psi$ in $\mathcal{D}_H(M)$. We can say that those solutions are $H$-weak solutions.
   However, in Lemma 6.1 in \cite{JFA} it is observed that  $H$-weak solutions to Equation (\ref{eqExample}) are weak solutions, that is $(\ref{weak-sol})$ also holds for every $\psi$ in $\mathcal{D}(M)$.

Therefore the solution given by
Theorem~\ref{reac-diff} is an $ H$-invariant weak solution in the usual unrestricted sense.

\begin{corollary}
\label{cor:unrestricted-weak-formulation}
Under the hypotheses of Theorem~\ref{reac-diff}, the $H$-invariant
solution $u\in L^{q(\cdot)}_{1,H}(M)$ obtained there is an $ H $-invariant weak solution of
\[
-\Delta_{q(x)}u+|u|^{q(x)-2}u=f(x,u)
\]
in the usual sense, namely
\[
 \int_M
 \left(
 |\nabla_g u|^{q(x)-2} g(\nabla_g u,\nabla_g \psi)  + |u|^{q(x)-2}u\psi  - f(x,u)\psi
 \right)
 dV_g  =0
\]
for every $\psi\in\mathcal D(M)$.
\end{corollary}

\subsection{Families of manifolds satisfying conditions of Theorem \ref{reac-diff}} \label{models}

\smallskip

Consider the warped product $n$-manifold $M = X \times Y$ equipped with a metric (\ref{warped}). When $X$ is 1-dimensional, the situation simplifies substantially, and expressions become 
easier to handle. We will use $(X,g_X) \equiv (\mathcal{R} , dr^2)$ as  the horizontal space, or base, with coordinate $r$, where $ \mathcal{R} $ is a complete and non-compact subset of $ \mathbb{R} $.
We recall once again that our manifolds may have non-empty boundary.

Let us assume that the function $q$ satisfies the hypotheses of Theorem \ref{reac-diff}, and that $u$ is an $H$-invariant solution of 
Equation (\ref{eqExample}).  Since $H$ acts transitively on $Y$, every $H$-invariant function on $M=\mathcal R\times Y$ depends only on $r$, and the same applies to
the exponents $q,p$, and to $f(\cdot,t)$ for each fixed $t$. Understanding that  the functions depend on the variable $r$, Equation (\ref{eqExample}) becomes the ODE

\begin{eqnarray}        \label{eqExample-1}
- \Delta_{\phi,q} u + \lvert  u \rvert^{q-2} u = f(r,u) \ ,
\end{eqnarray}
with

\begin{eqnarray}       \label{q-Laplacian}
\Delta_{\phi,q} u :=  \frac{1}{\phi^{n-1}}( \lvert  u_r \rvert^{q-2}  \phi^{n-1}\  u_r )_r  = \Delta_q u +
  (n-1) (\log \phi)_r  \lvert    u_r \rvert^{q-2} u_r   \ ,
 \end{eqnarray} 
 where 
$(\cdot)_r$ denotes the derivative of $(\cdot)$ with respect to $r$. 
In addition,  Equation (\ref{weak-sol}) becomes 
\begin{eqnarray}
\int_{\mathcal{R}}\  \Big(\   \lvert u_r \rvert^{q-2}  u_r \psi_r +   \lvert u \rvert^{q-2}  u \psi - f(r,u) \psi    \Big)
\phi^{n-1}\  dr = 0 \ ,
\end{eqnarray}
whenever $\psi$ is in $\mathcal{D}(\mathcal{R})$.

In Equation (\ref{q-Laplacian}),  we can formally write, if $q$ is $C^1$, for smooth functions $u$ with $u_r\neq 0$: 
 $$   \Delta_q  u = ( \lvert u_r \rvert^{q-2} u_r)_r  =  \lvert u_r \rvert^{q-2} u_{rr} +
   (\lvert u_r \rvert^{q-2})_r u_r \ , $$ 
   with 
    $$ ( \lvert u_r \rvert^{q-2})_r    =   \lvert u_r \rvert^{q-2} (\ q_r \log \lvert u_r \rvert + 
    (q-2) ( \log \lvert u_r \rvert)_r \ ) \ .$$
    Therefore 
$$ \Delta_{\phi,q} u  =  \lvert u_r \rvert^{q-2} \{\   u_{rr}  +    
 \{(\  (n-1) (\log \phi)_r  + q_r \log \lvert u_r \rvert + (q-2) ( \log \lvert u_r \rvert)_r  \  \} \ u_r \    \}  \ .  $$
 
\smallskip

We note that if $q \equiv 2$ is constant, and we consider $ \mathcal{R} \subset \mathbb{R}_+$, the choice  $\phi(r) = r$ allow us to recover the classic  radial Laplacian 
$$ \Delta u = u_{rr} + \frac{n-1}{r} u_r  $$
for the standard metric in $\mathbb{R}^n$ using spherical coordinates, and Equation (\ref{eqExample-1}) becomes the non-linear ODE
$$
- \left( u_{rr} + \frac{n-1}{r} u_r \right) + u = f(r,u)\; .
$$

\bigskip

In order to apply Theorem \ref{reac-diff}, $ \mathcal{R} \subseteq \mathbb{R} $ and $\phi : \mathcal{R} \to (0,\infty)$ must be such that:
\begin{itemize}
\item  Hypothesis \ref{bounded-geo} holds, i.e. $(M,g_{\phi})$ has property $\mathcal{B}_{\text{vol}}(\lambda,\vartheta)$ for some  $(\lambda, \vartheta )$. 
\item  Hypothesis \ref{M-balls} is satisfied. 
\end{itemize}

In the warped product examples considered below, hypothesis
\ref{M-balls} will be verified by showing that the $H$-orbits
$\{r\}\times Y$ contain arbitrarily many disjoint balls as
$r \to\infty$. 

\bigskip

Our first aim is to give sufficient conditions on $\phi$ to ensure that the lower bound for the Ricci tensor required  by hypothesis \ref{bounded-geo} in Theorem \ref{reac-diff} holds.

\begin{theorem}   \label{Warp-1D}
Consider the complete non-compact warped product $(M = \mathcal{R} \times Y , g_{\phi})$, with $g_{\phi} = dr^2  \oplus \phi^2 g_{Y}$  as in (\ref{warped}), where
 $(Y , g_Y)$ is an Einstein homogeneous space of dimension $n-1$ for the compact Lie group $H$, $g_Y$  being the Einstein metric such that  
 $  Rc(g_Y) = \lambda_Y  (n-2)   g_Y$, in which we assume $n\geq 3$ and that $\dim(Y) = n-1$. Then
\begin{equation}    \label{Rc-bound}
 Rc(g_{\phi})  \geq \lambda_M (n-1)    g_{\phi}  
\end{equation}
for a given constant number $\lambda_M \leq 0 \leq \lambda_Y $ whenever we have 
\begin{equation}   \label{Rc-hor}
 \lambda_M  \leq - \frac{\phi_{rr}}{\phi}
 \end{equation}
and 
\begin{equation}    \label{Rc-vert}
 - \phi^2 (\log \phi)_{rr} \leq \lambda_Y   \ .  
 \end{equation}

\end{theorem}

\smallskip

\begin{proof}

First  note that we do need the hypothesis $\lambda_M \leq 0 \leq \lambda_Y $: Indeed, if  $\lambda_M > 0$ (and $(M,g_\phi)$ is complete) then $M$ must be  compact;  also, $0 \leq \lambda_Y $ 
follows from Remark \ref{H-spaces}. We split (\ref{Rc-bound}) in  horizontal and vertical components, so using Lemma \ref{Ricci-warp}, we see  that (\ref{Rc-bound}) is equivalent to
$$  \lambda_M     \leq \ -  \frac{\phi_{rr}}{\phi}  $$
together with
\begin{eqnarray}     \label{Rc-vert-1}
\lambda_M \leq 
\frac{1}{n-1} \{\  \frac{\lambda_Y (n-2)}{ \phi^2} - \frac{1}{n-1} \frac{(\phi^{n-1})_{rr}}{\phi^{n-1}}\  \} \ ,
\end{eqnarray}
and we rewrite (\ref{Rc-vert-1})  as
\begin{eqnarray}   \label{Rc-vert-2} 
\frac{\phi^2}{n-2} \{\ (n-1) \lambda_M  + \frac{1}{n-1} \frac{(\phi^{n-1})_{rr}}{\phi^{n-1}}\  \} \leq \lambda_Y\ .
\end{eqnarray}

\noindent Using the identity 
$$  - (n-2) ( \log \phi )_{rr} = - (n-1) \frac{\phi_{rr}}{\phi} + \frac{1}{n-1} \frac{(\phi^{n-1})_{rr}}{\phi^{n-1}} \ , $$
from the inequality
$$ - \phi^2 (\log \phi)_{rr} \leq \lambda_Y $$
together with (\ref{Rc-hor}), we recover (\ref{Rc-vert-2}).
\end{proof}

The next corollary provides upper bounds for warping functions satisfying
the curvature inequalities of Theorem~\ref{Warp-1D}. The warp functions obtained in the case of equality
will be those of some \textit{standard} model spaces, see Remark \ref{standard} below.

\begin{corollary}   \label{coro-warp}
In the context of Theorem  $\ref{Warp-1D}$, assume that the warp function $\phi : \mathcal{R} \subset \mathbb{R}_+ \to (0,\infty)$ 

satisfies the curvature inequalities
\[
\lambda_M\leq -\frac{\phi_{rr}}{\phi}
\]
and
\[
-\phi^2(\log\phi)_{rr}\leq \lambda_Y .
\]
Then $\phi$ has the following upper bounds, in which we assume $r \geq 0$:

\begin{enumerate}
\item If $\lambda_M = 0 < \lambda_Y$, there is a constant $D > 0$ such that $ \phi(r) \leq C r + D$, where 
$C = \sqrt{\lambda_Y} > 0$.

\item If $\lambda_M < 0 = \lambda_Y $, then $\phi(r) \leq C e^{\sqrt{-\lambda_M}  r}$, where $C > 0$ is a constant.

\item If $\lambda_M < 0 < \lambda_Y$, there is a constant $D$ such that
$$ \phi(r) \leq  C \sinh \left(  \sqrt{-\lambda_M} r + D    \right) \ , $$
where $\sqrt{\frac{\lambda_Y}{-\lambda_M}} =  C  > 0$.

\item If $\lambda_M = 0 = \lambda_Y$, then there is a constant $ C > 0 $ such that
$$ \phi(r) =  C . $$ 

\end{enumerate}

\end{corollary}


\begin{proof}
Note that from (\ref{Rc-hor}) and (\ref{Rc-vert}) we get
\begin{equation} \label{Rc-hor-ver}
\lambda_M \phi^2 + (\phi_r)^2 \leq \lambda_Y \ .
\end{equation}
\begin{enumerate}
\item If $\lambda_M = 0$, integrating (\ref{Rc-hor-ver}) we have $\phi(r) \leq C r + D$, where $C = \sqrt{\lambda_Y}$.

\item If $\lambda_Y = 0$, then (\ref{Rc-hor-ver}) gives $\phi(r) \leq C e^{\sqrt{-\lambda_M}  r}$, where $C$ is a constant.

\item If $\lambda_M < 0 < \lambda_Y$, from (\ref{Rc-hor-ver}) we obtain
$  |\phi_r| / \sqrt{\lambda_Y - \lambda_M \phi^2} \leq 1 $, so that
\[
\left|
\frac{d}{dr}
\sinh^{-1}\left(
\sqrt{\frac{-\lambda_M}{\lambda_Y}}\phi(r)
\right)
\right|
\leq
\sqrt{-\lambda_M}.
\]
Integrating this inequality on $[r_0,r]$ gives the claimed bound.

\item If $\lambda_M = 0 = \lambda_Y$, then (\ref{Rc-hor-ver}) gives $ \phi(r)_r = 0 $.
\end{enumerate}
\end{proof}  


\begin{remark}  \label{not-cosh}
A case excluded in  Theorem \ref{Warp-1D} is when $\lambda_M < 0 $ and  $\lambda_Y < 0 $, as in this case we would obtain a non-compact fiber $Y$, see Remark \ref{H-spaces}. 
In this case, using (\ref{Rc-hor-ver}) we have 
\[ \left| \frac{d}{dr} \cosh^{-1}\left(
\sqrt{\frac{\lambda_M}{\lambda_Y}}\phi(r)  \right) \right| = 
\frac{ | \sqrt{\frac{\lambda_M}{\lambda_Y}}. \phi_r| }{   \sqrt{ -1  + \frac{\lambda_M}{\lambda_Y} \phi^2}} \leq \sqrt{-\lambda_M} ,\]
and integrating between $ [r_0 , r ] $ we obtain
\[   \phi(r) \leq  C \cosh \left(  \sqrt{-\lambda_M} r + D    \right) , \]
where $\sqrt{\frac{\lambda_Y}{\lambda_M}} =  C  > 0$. 
\end{remark}


\begin{remark}     \label{standard}
When equality holds in the curvature inequalities of Theorem~\ref{Warp-1D} and Remark \ref{not-cosh},
the model warping functions are precisely the extremal functions appearing
in Corollary~\ref{coro-warp}. This gives the Einstein model
spaces: compare with 9.109 in \cite{Besse}, and Table 1 in \cite{Solitons}.
\end{remark}

To verify that hypothesis \ref{bounded-geo}  in Theorem \ref{reac-diff} holds, in addition to the conditions of  Theorem \ref{Warp-1D}, we must have the lower bound 
$\lvert B_1(x) \rvert_{g_{\phi}} \geq \vartheta $ for the volume of unit balls in M. For these purposes, recall that
the Riemannian measure of $g_\phi$ is
\[
dV_{g_\phi}(r,y)=\phi(r)^{n-1}\ dr \wedge dV_{g_Y}(y),
\]
and that the distance between points on $ Y $ in the fiber at $ r \in \mathcal{R} $ scales as $  \phi(r)\   d_Y $, where $ d_Y $ is the distance in $ (Y, g_Y) $.\\

We will assume that the complete subset $ \mathcal{R} \subset \mathbb{R} $ is such that:
\begin{enumerate}
\item \label{A} $ \inf\ \{\ \phi(r) : r \in \mathcal{R}\ \} \geq  1/(\ 2\ \text{inj} (Y, g_Y)\ )  > 0      $, where $ \text{inj} (Y, g_Y)$ denotes the injectivity radius of $ (Y,g_Y) $, and
\item \label{B}  There exists a constant   $ L_{\phi} < \infty $  such that  $ \sup\ \{\ \lvert ( \log \phi(r) )_r \rvert  : r \in \mathcal{R}\ \}  \leq L_{\phi} $.  
\end{enumerate}
Fix a point $ x= (r , y)  $ in $ M $, and consider the set
\[ Q_{\frac{1}{2}}(r,y) := \{  (s,z) \in \mathcal{R} \times Y :  \vert  s- r   \vert < \frac{1}{2}\  ,\  \phi(r)\ d_Y(y,z) < \frac{1}{2}        \}  = \{ s \in \mathcal{R} :  \vert  s- r   \vert < \frac{1}{2} \} \times  B^{d_Y}_{1 / 2 \phi(r) }(y) ,  \]
where $ B^{d_Y}_R(y) $ is the ball of radius $ R $ in $ (Y , g_Y) $ centered at y.\footnote{Note that 
$ B^{\phi(r) d_{Y}}_R(y) = B^{d_Y}_{R / \phi(r) }(y) $ .}\\

Therefore $ Q_{\frac{1}{2}}(r,y) $ is contained in $ B_1(r,y) $, hence, denoting here and henceforth the $g$-volume of a ball $ B $ by $ \lvert B \rvert_g $,
\[ \lvert B_1(r,y) \rvert_{g_{\phi}}  \geq \lvert  Q_{\frac{1}{2}}(r,y) \rvert_{g_{\phi}}   = \int_{Q_{\frac{1}{2}}(r,y) } d V_{g_{\phi}} \geq   \int_{r}^{r+\frac{1}{2}} \phi(s)^{n-1} ds\ \lvert  B^{d_Y}_{1 / 2 \phi(r) } \rvert_{g_Y} ,\  \]
where the last inequality is necessary in order to take into account  the case that $ r $ is at a distance less than $ 1/2 $ from the boundary of $ \mathcal{R} $, and the center of the ball  has been omitted by homogeneity of the volume in $ ( Y , g_Y ) $.

From assumption \ref{B}  we obtain that $ \phi(s) \geq e^{-{L_{\phi}/2}} \phi(r) $, therefore
\[ \lvert B_1(r,y) \rvert_{g_{\phi}}  \geq \frac{1}{2}\  e^{\frac{(1-n) L_{\phi}}{2}} \phi(r)^{n-1}\ \lvert B^{d_Y}_{1/2 \phi(r) } \rvert_{g_Y}\   .     \]

But \ref{A} says that $ 1 / 2 \phi(r) \leq \text{inj}(Y,g_Y) $, thus by Bishop-Gromov comparison we get \[ \lvert B^{d_Y}_{1 /2 \phi(r) } \rvert_{g_Y}   \geq  \omega_Y ( 1 / 2 \phi(r)  )^{n-1}   \] with 
 $ \omega_Y := \lvert   B^{d_Y}_{\text{inj}(Y, g_Y)} \rvert_{g_Y} / (\text{inj}(Y, g_Y))^{n-1}   > 0$,
to obtain 
\[   \lvert B_1(x) \rvert_{g_{\phi}} \geq    e^{\frac{(1-n) L_{\phi}}{2}} \omega_Y / 2^{n}            =: \vartheta   ,  \]
and the lower bound for the volume of balls in hypothesis 2 is verified. \\

Now we give sufficient conditions on $ \phi $ for hypothesis 3 to hold. Since $ H $ acts transitively on $ Y $, every subset of $ Y $ is contained in the orbit of $ H $. Our aim is, roughly speaking, to give conditions on $ \phi(r) $ to ensure that  there is an increasing number of disjoint balls of radius $ R $  in $ ( Y , \phi(r)^2 g_Y) $ as $  r  $ increases, such a  number going to infinity as  $  r  \to \infty  $. 
For those purposes, let $ \text{Pack}_{d_Y}(R) $ denote the maximal number of disjoint balls of radius $ R $ that can be packed in $ ( Y ,g_Y ) $; recalling that given $ R > 0,\  r \in \mathcal{R}$ and $ y \in Y $ one has that $ B_{R}^{\phi(r) d_Y} (y) = B_{R/ \phi(r)}^{d_Y}(y) $ as subsets of $ Y $, it follows that 
\[ 
 \text{Pack}_{\phi(r) d_Y}(R) = \text{Pack}_{d_Y}(R / \phi(r))       
 \]
is the maximal number of disjoint balls of radius $ R $ in $ ( Y ,\phi(r)^2 g_Y ) $.
To obtain a lower bound for $ \text{Pack}_{d_Y}(R) $, denote by $ \text{Cov}_{d_Y}(R) $
the minimal number of balls of radius $R $ needed to cover $ Y $.
It is not difficult to see that
\[
 \text{Cov}_{d_Y}(2 R)  \leq \text{Pack}_{d_Y}(R) \leq \text{Cov}_{d_Y}(R)  . 
 \]

Since $ (Y, g_Y ) $ is homogeneous of non-negative curvature, then
\[    
  \lvert Y \rvert_{g_Y} \leq \text{Cov}_{d_Y}(R)\ \lvert B^{d_Y}_R \rvert_{g_Y}  , \]
  and we have
  \[  \text{Pack}_{d_Y}(R)  \geq   \frac{\lvert Y \rvert_{g_Y} }{\lvert B^{d_Y}_{2R} \rvert_{g_Y} }  .\]
  
  However, since $ \text{Rc}(g_{Y}) \geq \lambda_Y (n-2) g_Y ,$ then by Bishop's estimate
  \[    \lvert B^{d_Y}_{R} \rvert_{g_Y} \leq   c(n) \int_0^R \Big( \frac{\sin (t \sqrt{\lambda_Y} )}{\sqrt{\lambda_Y}} \Big)^{n-2}\ dt \leq c(n) \int_0^R t^{n-2} dt = c'(n) R^{n-1} ,   \]
  compiling
  \[  \text{Pack}_{d_Y}(R)  \geq   \frac{\lvert Y \rvert_{g_Y} }{c'(n) 2^{n-1}} R^{1-n} =: C_Y R^{1-n}  ,\]
 and the estimate 
\[
\text{Pack}_{\phi(r) d_Y}(R) \geq C_Y R^{1-n} \phi(r)^{n-1}
\]
is obtained. Therefore,  the orbit of $ H $ has at least    
$    C_Y R^{1-n} \phi(r)^{n-1}      $ disjoint balls of radius $ R $ in the strip $ [ r - R , r + R ] \times Y $ inside $ \mathcal{R} \times Y $. \\

We infer that if $ \phi(r) \to \infty $ when $ r \to \infty $, then hypothesis 3 is satisfied. It will be convenient to summarize what has been shown.


\begin{theorem}
 Assume that $ (Y , g_Y) $ is a $ (n-1)$-dimensional compact manifold that  is  a homogeneous space for a Lie group $ H $, with $ Rc(g_Y) = \lambda_Y (n-2) g_Y $, $ \lambda_Y \geq 0 $. Consider a non-compact complete subinterval $ \mathcal{R} $  of $ \mathbb{R}_+ $, together with  a function $ \phi : \mathcal{R} \to \mathbb{R}_+ $ such that:
\begin{enumerate}
\item $ \lambda_M  \leq - \frac{\phi_{rr}}{\phi} $ , where $ \lambda_M \leq 0 $.
\item $ - \phi^2 (\log \phi)_{rr} \leq \lambda_Y   $. 
\item $ \inf\ \{\ \phi(r) : r \in \mathcal{R}\ \} \geq  1/(\ 2\ \text{inj} (Y, g_Y)\ )  > 0      $, where $ \text{inj} (Y, g_Y)$ denotes the injectivity radius of $ (Y,g_Y) $.
\item  $ \sup\ \{\ \lvert ( \log \phi(r) )_r \rvert  : r \in \mathcal{R}\ \}  \leq L_{\phi} $ for a constant   $ L_{\phi} < \infty $.
\item $ \phi(r) \to \infty $ when $ r \to \infty $.  
\end{enumerate}

Then  $ ( M , g_{\phi} ) = (\mathcal{R} \times Y , dr^2 \oplus \phi^2 g_Y ) $ satisfies hypotheses 2 and 3 of Theorem 1, with $ Rc(g_{\phi}) \geq \lambda_M (n-1) g_{\phi}. $ 
\end{theorem}

\begin{remark}
All the standard model spaces  in  Corollary \ref{coro-warp} satisfy conditions 1, 2, 3, 4 in a suitable interval $ \mathcal{R} = [ r_0 , \infty ) \subset \mathbb{R}_+ $. Only in the case when 
$ \lambda_M = 0 = \lambda_Y $
condition 5 is not satisfied, for then $ \phi(r) $ is constant.\\
\end{remark}






\section{Radial reductions and model solutions}\label{sec:radial}

Throughout this section we write $m:=\dim Y=n-1$, so that $\dim M=m+1$.
Here $k$ denotes the constant diffusion exponent and is independent of the
fibre dimension.

Let $\mathcal{R} \subset\mathbb R$ be an interval. We consider radial equations on
\[
\mathcal{R} \times_\phi Y,
\qquad
g_\phi=dr^2+\phi(r)^2g_Y.
\]

We discuss explicit radial model solutions and logarithmic
reductions of the $H$-invariant equation (\ref{eqExample}) on $( M = \mathcal{R} \times Y , g_{\phi}  )$ obtained in Section~3. These
solutions should be understood primarily as strong solutions of the reduced
one-dimensional equation. Membership in the variational space
$L_1^{q(\cdot)}(M)$ must be checked separately, and in general the explicit
constant solutions below do not belong to $L_1^{q(\cdot)}(M)$ when
$(M,g_\phi)$ has infinite volume.

In full generality, the equation of interest for us is
\begin{eqnarray}        \label{eqExample-11}
- \Delta_{\phi,q} u + \lvert  u \rvert^{q-2} u = f(r,u) \ ,
\end{eqnarray}
in which $q$ and $f$ may depend on $r \in \mathcal{R}$, with
 
\begin{eqnarray}       \label{q-Laplacian-1}
\Delta_{\phi,q}u =
\left(|u_r|^{q(r)-2}u_r\right)_r
+
(n-1)(\log\phi)_r |u_r|^{q(r)-2}u_r\; .
\end{eqnarray}

\subsection{$q=2$}
From (\ref{q-Laplacian-1}) we obtain
$$
\Delta_{\phi,2} u  = \Delta u + (n-1) (\log \phi)_r  u_r
$$
and Equation (\ref{eqExample-11}) becomes
\begin{equation}        \label{eqExample-111}
- u_{rr} - (n-1) (\log \phi)_r  u_r +  u = f(r, u) \ .
\end{equation}

This is the classical radial semilinear equation associated with the warped
product measure $\phi(r)^{n-1}dr$. We shall not pursue this semilinear case
further, since the nonlinear logarithmic reduction below is the one relevant
for the $k$-Laplacian model.

\subsection{$q = k > 2$}

\begin{theorem}[Logarithmic reduction for the case $q=k>2$]
\label{thm:qeqk_log}
Let $n\ge 2$ and $k>2$. Assume $\phi\in C^{1}((0,\infty))$ with $\phi(r)>0$ for all $r>0$.
Let $\gamma,p:(0,\infty)\to\mathbb{R}$ be measurable functions.
Consider the radial equation
\begin{equation}
\label{eq:main_qeqk}
-\Delta_{\phi,k}u + u^{k-1} \;=\; \gamma(r)\,u^{p(r)-1}\qquad\text{a.e. on }(0,\infty),
\end{equation}
where for a radial function $u=u(r)$ we set
\[
\Delta_{\phi,k}u
:=\frac{1}{\phi^{n-1}(r)}\Big(\phi^{n-1}(r)\,|u_r|^{k-2}u_r\Big)_r
\quad\text{(weak derivative in $r$)}.
\]
Assume $u:(0,\infty)\to(0,\infty)$ satisfies
\[
u\in C^{1}((0,\infty)),
\qquad
\phi^{n-1}|u_r|^{k-2}u_r\in W^{1,1}_{\mathrm{loc}}((0,\infty)).
\]
Define
\[
v(r):=(k-1)\log u(r),
\qquad
\alpha(r):=\Big(\frac{v(r)}{k-1}\Big)_r=\frac{v_r(r)}{k-1},
\qquad
\beta(r):=|\alpha(r)|^{k-2}\alpha(r).
\]
Then \eqref{eq:main_qeqk} holds if and only if
\begin{equation}
\label{eq:beta_eq}
\beta_r+\Big((n-1)(\log\phi)_r+v_r\Big)\beta
=
1-\gamma(r)\exp\!\Big(\frac{p(r)-k}{k-1}\,v\Big)
\qquad\text{a.e. on }(0,\infty).
\end{equation}
Equivalently, in divergence form,
\begin{equation}
\label{eq:div_form}
\Big(\phi^{n-1}e^{v}\beta\Big)_r
=
\phi^{n-1}\Big(e^{v}-\gamma(r)\exp\!\big(\tfrac{p(r)-1}{k-1}\,v\big)\Big)
\qquad\text{a.e. on }(0,\infty).
\end{equation}
Moreover, since $\beta=(k-1)^{1-k}|v_r|^{k-2}v_r$, equation \eqref{eq:beta_eq} is equivalent to
\begin{equation}
\label{eq:v_second}
(|v_r|^{k-2}v_r)_r+\Big((n-1)(\log\phi)_r+v_r\Big)|v_r|^{k-2}v_r
=
(k-1)^{k-1}\Big[1-\gamma(r)\exp\!\Big(\frac{p(r)-k}{k-1}\,v\Big)\Big]
\end{equation}
in the sense of distributions on $(0,\infty)$ (and pointwise a.e.\ if $v$ is $C^2$).
\end{theorem}

\begin{proof}
We use only identities valid almost everywhere, with derivatives interpreted in the weak sense when needed.

\medskip
Firstly, introduce the logarithmic change of variables: Since $u>0$ and $u\in C^1$, the function $v=(k-1)\log u$ belongs to $C^1$ and
\begin{equation}
\label{eq:ur_identity}
v_r=(k-1)\frac{u_r}{u}
\qquad\Longleftrightarrow\qquad
u_r=\frac{v_r}{k-1}\,u=\alpha\,u.
\end{equation}
Also,
\begin{equation}
\label{eq:ukminus1}
u^{k-1}=e^{(k-1)\log u}=e^{v}.
\end{equation}
Using \eqref{eq:ur_identity} and $u>0$ we obtain
\begin{equation}
\label{eq:flux_identity}
|u_r|^{k-2}u_r
=
|\alpha u|^{k-2}(\alpha u)
=
|\alpha|^{k-2}\alpha\,u^{k-1}
=
\beta\,e^{v}.
\end{equation}

\medskip

By definition of $\Delta_{\phi,k}$,
\[
-\Delta_{\phi,k}u+u^{k-1}
=
-\frac{1}{\phi^{n-1}}\Big(\phi^{n-1}|u_r|^{k-2}u_r\Big)_r + u^{k-1}.
\]
Multiply \eqref{eq:main_qeqk} by $\phi^{n-1}$ and use \eqref{eq:ukminus1}--\eqref{eq:flux_identity} to obtain the equivalent identity
\begin{equation}
\label{eq:multiplied}
-\Big(\phi^{n-1}e^{v}\beta\Big)_r+\phi^{n-1}e^{v}
=
\phi^{n-1}\gamma(r)\,u^{p(r)-1}.
\end{equation}
Since $u=e^{v/(k-1)}$, we have
\begin{equation}
\label{eq:upminus1}
u^{p(r)-1}=\exp\!\Big(\frac{p(r)-1}{k-1}\,v\Big).
\end{equation}
Substituting \eqref{eq:upminus1} into \eqref{eq:multiplied} and rearranging yields \eqref{eq:div_form}.
Thus \eqref{eq:main_qeqk} implies \eqref{eq:div_form}.

Conversely, if \eqref{eq:div_form} holds and we define $u:=e^{v/(k-1)}$, then \eqref{eq:ukminus1}--\eqref{eq:flux_identity} hold identically, and reversing the above algebra shows \eqref{eq:multiplied}, hence \eqref{eq:main_qeqk}.
Therefore \eqref{eq:main_qeqk} and \eqref{eq:div_form} are equivalent.

\medskip
Now, expanding the derivative and obtaining the first-order equation for $\beta$: 
Set $\mu(r):=\phi^{n-1}(r)e^{v(r)}$. Then $\mu\in C^1((0,\infty))$ and $\mu>0$.
By hypothesis, $\phi^{n-1}|u_r|^{k-2}u_r\in W^{1,1}_{\mathrm{loc}}$; by \eqref{eq:flux_identity} this is exactly
$\mu\beta\in W^{1,1}_{\mathrm{loc}}((0,\infty))$.

We now use the standard product rule: if $a\in C^1_{\mathrm{loc}}$ and $b\in W^{1,1}_{\mathrm{loc}}$, then
$ab\in W^{1,1}_{\mathrm{loc}}$ and $(ab)_r=a_rb+ab_r$ a.e. We wish to apply this rule with $a=\mu$ and $b=\beta$ (note that $b$ is a.e. defined because 
$\mu\beta\in W^{1,1}_{\mathrm{loc}}$ and $\mu>0$). We justify the required Sobolev regularity of $\beta$.  By assumption,
\[
\mu\beta=\phi^{n-1}e^v\beta\in W^{1,1}_{\mathrm{loc}}(\mathcal{R}).
\]
Moreover, since $\phi>0$ and $v$ is locally absolutely continuous, we have
\[
\mu>0,
\qquad
\frac{1}{\mu}\in C^1_{\mathrm{loc}}(\mathcal{R})
\]
in the smooth case, and in particular $1/\mu\in W^{1,\infty}_{\mathrm{loc}}(\mathcal{R})$
under the corresponding local regularity assumptions. Hence
\[
\beta
=
\frac{\mu\beta}{\mu}
\in W^{1,1}_{\mathrm{loc}}(\mathcal{R}).
\]
Therefore the product rule may be applied in the weak sense, giving
\begin{equation}
\label{eq:product_rule_mu}
(\mu\beta)_r=\mu_r\beta+\mu\beta_r
\qquad\text{a.e. on }(0,\infty).
\end{equation}
Compute $\mu_r/\mu$:
\[
\frac{\mu_r}{\mu}
=
\frac{(\phi^{n-1})_r}{\phi^{n-1}}+\frac{(e^{v})_r}{e^{v}}
=
(n-1)(\log\phi)_r+v_r.
\]
Divide \eqref{eq:product_rule_mu} by $\mu>0$ and obtain
\begin{equation}
\label{eq:expand}
\frac{(\mu\beta)_r}{\mu}
=
\beta_r+\Big((n-1)(\log\phi)_r+v_r\Big)\beta
\qquad\text{a.e.}
\end{equation}
Now divide \eqref{eq:div_form} by $\mu=\phi^{n-1}e^{v}$ (which is strictly positive) to get
\[
\frac{(\mu\beta)_r}{\mu}
=
1-\gamma(r)\exp\!\Big(\frac{p(r)-1}{k-1}\,v-v\Big)
=
1-\gamma(r)\exp\!\Big(\frac{p(r)-k}{k-1}\,v\Big)
\qquad\text{a.e.}
\]
Together with \eqref{eq:expand}, this is exactly \eqref{eq:beta_eq}.
Thus \eqref{eq:main_qeqk}$\iff$ \eqref{eq:div_form} $\iff$ \eqref{eq:beta_eq}.

\medskip
Now, we show the equivalence with the second-order equation for $v$.
By definition $\alpha=v_r/(k-1)$, hence
\[
\beta=|\alpha|^{k-2}\alpha
=
\left|\frac{v_r}{k-1}\right|^{k-2}\frac{v_r}{k-1}
=
(k-1)^{1-k}|v_r|^{k-2}v_r
\qquad\text{a.e.}
\]
Therefore $\beta_r=(k-1)^{1-k}(|v_r|^{k-2}v_r)_r$ in the sense of distributions (and a.e.\ if $v\in C^2$).
Substituting these identities into \eqref{eq:beta_eq} and multiplying by $(k-1)^{k-1}$ yields \eqref{eq:v_second}.
This completes the proof.
\end{proof}

\begin{example}[Constant strong solutions of the reduced equation]
Assume, in addition to the hypotheses of Theorem~\ref{thm:qeqk_log}, that
\[
p(r)\equiv p_0\in\mathbb{R},\qquad \gamma(r)\equiv \gamma_0\in\mathbb{R}.
\]
Then the following give explicit (radial) solutions of \eqref{eq:main_qeqk}:

\begin{enumerate}
\item If $p_0\neq k$ and $\gamma_0>0$, then
\begin{equation}
\label{eq:const_u_pk}
u(r)\equiv u_0:=\gamma_0^{\frac{1}{\,k-p_0\,}}>0
\end{equation}
is a solution.

\item If $p_0=k$ and $\gamma_0=1$, then for every $C>0$,
\begin{equation}
\label{eq:const_u_pk_eq}
u(r)\equiv C
\end{equation}
is a solution.
\end{enumerate}
In both cases, the associated logarithmic variables satisfy
\[
v(r)=(k-1)\log u(r)\equiv \text{constant},\qquad \alpha(r)=\frac{v_r}{k-1}\equiv 0,\qquad \beta(r)=|\alpha|^{k-2}\alpha\equiv 0.
\]
\end{example}

\begin{proof}
Let $u$ be any positive constant. Then $u_r\equiv 0$, hence $|u_r|^{k-2}u_r\equiv 0$ and therefore
\[
\Delta_{\phi,k}u=\frac{1}{\phi^{n-1}}\big(\phi^{n-1}\cdot 0\big)_r\equiv 0
\quad\text{on }(0,\infty).
\]
Thus \eqref{eq:main_qeqk} reduces to the algebraic identity
\[
u^{k-1}=\gamma_0\,u^{p_0-1}.
\]
If $p_0\neq k$, this is equivalent (since $u>0$) to $u^{k-p_0}=\gamma_0$. A positive solution exists exactly when $\gamma_0>0$, and is uniquely given by \eqref{eq:const_u_pk}.
If $p_0=k$, the identity becomes $u^{k-1}=\gamma_0 u^{k-1}$, which holds for all $u>0$ if and only if $\gamma_0=1$, yielding \eqref{eq:const_u_pk_eq}.

Finally, for any positive constant $u$, one has $v=(k-1)\log u\equiv\text{const}$, hence $v_r\equiv 0$ and consequently $\alpha\equiv 0$ and $\beta\equiv 0$.
\end{proof}

\begin{remark}
If the complete warped product has infinite volume, as it happens if we use as warping functions the upper bounds appearing in cases 1, 2, 3 of Corollary 3.3,
the positive constant solutions above do
not belong to $L^k(M)$. Hence they should not be interpreted as variational
solutions given by Theorem~\ref{reac-diff}; they are only strong solutions
of the reduced equation.
\end{remark}

\begin{proposition}
\label{prop:nonconstant_family}
Assume the hypotheses of Theorem~\ref{thm:qeqk_log} and, in addition,
\[
p(r)\equiv k,\qquad \gamma(r)\equiv 1 \qquad \text{on }(0,\infty).
\]
Then \eqref{eq:main_qeqk} reduces to the $k$--harmonic radial equation
\begin{equation}
\label{eq:kharmonic}
\big(\phi^{n-1}(r)\,|u_r|^{k-2}u_r\big)_r=0\qquad \text{on }(0,\infty)
\end{equation}
(in the distributional sense, hence a.e.).
Consequently, for any $r_0>0$ and any constants $A\in\mathbb R$, $B\in\mathbb R$, the function
\begin{equation}
\label{eq:general_phi_solution}
u(r):=A+B\int_{r_0}^{r}\phi(s)^{-\frac{n-1}{k-1}}\,ds
\end{equation}
satisfies \eqref{eq:kharmonic} on $(0,\infty)$, and it is non-constant whenever $B\neq 0$.
Moreover, if $\phi(r)=r$ (Euclidean case), then \eqref{eq:general_phi_solution} becomes
\begin{equation}
\label{eq:euclidean_solution}
u(r)=
\begin{cases}
A+\displaystyle B\,\frac{k-1}{k-n}\,r^{\frac{k-n}{k-1}}, & k\neq n,\\[1.2ex]
A+B\log r, & k=n,
\end{cases}
\end{equation}
and any choice of $A,B$ for which $u(r)>0$ on the interval of interest yields a positive non-constant solution of \eqref{eq:main_qeqk}.
\end{proposition}

\begin{proof}
Assume first that $u$ is a positive radial solution in the resonant case
\[
p\equiv k,
\qquad
\gamma\equiv 1.
\]
Then the equation
\[
-\Delta_{\phi,k}u+u^{k-1}=u^{k-1}
\]
reduces to the radial $k$-harmonic equation
\[
\Delta_{\phi,k}u=0.
\]
Using the radial formula
\[
\Delta_{\phi,k}u
=
\frac{1}{\phi^{n-1}}
\left(
\phi^{n-1}|u_r|^{k-2}u_r
\right)_r,
\]
we obtain
\[
\left(
\phi^{n-1}|u_r|^{k-2}u_r
\right)_r=0.
\]
Hence there exists a constant $C\in\mathbb R$ such that
\[
\phi^{n-1}|u_r|^{k-2}u_r=C.
\]
Equivalently,
\[
|u_r|^{k-2}u_r
=
C\,\phi^{-(n-1)}.
\]

Now recall that the map
\[
\mathbb R\ni s\longmapsto |s|^{k-2}s\in\mathbb R
\]
is a strictly increasing bijection for every $k>1$. Therefore there exists
a unique constant $B\in\mathbb R$ such that
\[
|B|^{k-2}B=C.
\]
Explicitly,
\[
B=
\begin{cases}
\operatorname{sgn}(C)|C|^{1/(k-1)}, & C\neq0,\\
0, & C=0.
\end{cases}
\]
It follows that
\[
u_r
=
B\,\phi^{-\frac{n-1}{k-1}}.
\]
Integrating from a fixed point $r_0\in I$, we obtain
\[
u(r)
=
A+
B\int_{r_0}^r
\phi(s)^{-\frac{n-1}{k-1}}\,ds,
\]
where $A=u(r_0)$.

Conversely, let
\[
u(r)
=
A+
B\int_{r_0}^r
\phi(s)^{-\frac{n-1}{k-1}}\,ds
\]
and assume that $u>0$ on the interval under consideration. Then
\[
u_r
=
B\,\phi^{-\frac{n-1}{k-1}}.
\]
Consequently,
\[
|u_r|^{k-2}u_r
=
|B|^{k-2}B\,
\phi^{-\frac{(n-1)(k-1)}{k-1}}
=
|B|^{k-2}B\,\phi^{-(n-1)}.
\]
Multiplying by $\phi^{n-1}$, we get
\[
\phi^{n-1}|u_r|^{k-2}u_r
=
|B|^{k-2}B,
\]
which is constant. Hence
\[
\left(
\phi^{n-1}|u_r|^{k-2}u_r
\right)_r=0.
\]
Therefore
\[
\Delta_{\phi,k}u=0.
\]
Since $p=k$ and $\gamma=1$, this is equivalent to
\[
-\Delta_{\phi,k}u+u^{k-1}=u^{k-1}
=
\gamma u^{p-1}.
\]
Thus $u$ solves the resonant radial equation. The proof is complete.
\end{proof}

The family in Proposition~\ref{prop:nonconstant_family} gives strong radial
solutions of the reduced equation, but regretfully, on complete warped products of infinite
volume, these solutions need not belong to $L_1^k(M)$, as we observe in the following proposition.

\begin{proposition}[Non-integrability of the constant radial profiles]
Let $Y$ be compact and let $a>0$. If either
\[
M=\mathcal{R}\times_{e^{ar}}Y,
\]
then every positive constant function $u\equiv C>0$ satisfies
\[
u\notin L^q(M)
\]
for every $q>0$. In particular, the constant strong solutions constructed above do not belong to $L_1^q(M)$.
\end{proposition}
\begin{proof}
Let $m = n-1 =\dim Y$. The Riemannian measure is
\[
dV_{g_\phi}=\phi(r)^m\,dr \wedge dV_{g_Y}.
\]
Since $Y$ has positive finite volume, we have
\[
\operatorname{Vol}_{g_\phi}(M)
=
\operatorname{Vol}_{g_Y}(Y)
\int_{\mathcal{R}} \phi(r)^m\,dr,
\]
with $ \mathcal{R}=[r_0,\infty), r_0 > 0$.

If $\phi(r)=e^{ar}$, then
\[
\int_{\mathcal{R}}\phi(r)^m\,dr
=
\int_{\mathcal{R}} e^{amr}\,dr
=
\infty.
\]
If $\phi(r)=ar$, then
\[
\int_{\mathcal{R}} \phi(r)^m\ dr
=
a^m\int_{\mathcal{R}} r^m\,dr
=
\infty.
\]
Thus, in both cases,
\[
\operatorname{Vol}_{g_\phi}(M)=\infty.
\]
Therefore, for $u\equiv C>0$,
\[
\int_M |u|^q\,dV_{g_\phi}
=
C^q\operatorname{Vol}_{g_\phi}(M)
=
\infty.
\]
Hence $u\notin L^q(M)$ for every $q>0$. In particular, for every
$q>1$, $u\notin L_1^q(M)$.

\end{proof}

\subsection{The case $q=q(r)$}

We now consider the case in which the exponent $q$ depends on the radial
variable. The calculation is analogous to the constant exponent case, but the
logarithmic change of variables must be adapted to the pointwise exponent
$q(r)$. We present the result in a smooth setting, since the purpose of
this subsection is to describe the reduced equation satisfied by the
logarithmic variable.

\begin{proposition}
\label{prop:q-variable-log-reduction}
Let $ \mathcal{R} \subset\mathbb R$ be an interval as before, let $\phi\in C^1(\mathcal{R}), \phi(r)>0$ and let $q,p\in C^1(\mathcal{R}), q(r)>1$. 
Let $\gamma: \mathcal{R}\to\mathbb R$ be a measurable function. Consider the radial
equation
\begin{equation}
\label{eq:main_q_variable}
-\Delta_{\phi,q}u+u^{q(r)-1}
=
\gamma(r)u^{p(r)-1},
\end{equation}
where
\begin{equation}
\label{eq:phi_q_laplacian_variable}
\Delta_{\phi,q}u
:=
\frac{1}{\phi(r)^{n-1}}
\left(
\phi(r)^{n-1}|u_r|^{q(r)-2}u_r
\right)_r .
\end{equation}
Assume that $u: \mathcal{R} \to(0,\infty)$ is smooth. Define the logarithmic variable
$v$ by
\begin{equation}
\label{eq:q_variable_log_change}
u(r)=\exp\left(\frac{v(r)}{q(r)-1}\right)\; ,
\end{equation}
or, equivalently, $v(r)=(q(r)-1)\log u(r)$. 
Define
\begin{equation}
\label{eq:alpha_beta_variable}
\alpha(r):=
\left(\frac{v(r)}{q(r)-1}\right)_r,
\qquad
\beta(r):=|\alpha(r)|^{q(r)-2}\alpha(r).
\end{equation}
Then, if in addition $\beta\in C^1(\mathcal{R})$, $u$ solves \eqref{eq:main_q_variable} if and only if
\begin{equation}
\label{q-ODE}
\beta_r
+
\beta\left\{(n-1)(\log\phi)_r+v_r\right\}
=
1-\Gamma(r),
\end{equation}
where
\begin{equation}
\label{eq:Gamma_variable}
\Gamma(r):=
\gamma(r)
\exp\left(
\frac{p(r)-q(r)}{q(r)-1}v(r)
\right).
\end{equation}
Equivalently, \eqref{q-ODE} may be written as
\begin{equation}
\label{eq:q_variable_beta_equation}
\beta_r
+
\left((n-1)(\log\phi)_r+v_r\right)\beta
=
1-
\gamma(r)
\exp\left(
\frac{p(r)-q(r)}{q(r)-1}v(r)
\right)\; ,
\end{equation}
or, in divergence form, 
\begin{equation}
\label{eq:q_variable_divergence_form}
\left(
\phi^{n-1}e^v\beta
\right)_r
=
\phi^{n-1}
\left[
e^v
-
\gamma(r)
\exp\left(
\frac{p(r)-1}{q(r)-1}v
\right)
\right].
\end{equation}
\end{proposition}

\begin{proof}
We work directly with the reduced radial equation. For radial functions,
\[
\Delta_{\phi,q}u
=
\frac{1}{\phi^{n-1}}
\left(
\phi^{n-1}|u_r|^{q(r)-2}u_r
\right)_r .
\]
Equivalently, in expanded form,
\[
\Delta_{\phi,q}u
=
\left(|u_r|^{q(r)-2}u_r\right)_r
+
(n-1)(\log\phi)_r |u_r|^{q(r)-2}u_r .
\]

Let
\[
u(r)=e^{\frac{v(r)}{q(r)-1}}.
\]
Then
\[
u_r
=
\left(\frac{v}{q-1}\right)_r u
=
\alpha u,
\]
where
\[
\alpha(r)=\left(\frac{v(r)}{q(r)-1}\right)_r.
\]
Moreover,
\[
u^{q(r)-1}=e^v.
\]
Hence
\[
|u_r|^{q(r)-2}u_r
=
|\alpha u|^{q(r)-2}\alpha u
=
|\alpha|^{q(r)-2}\alpha\,u^{q(r)-1}
=
\beta e^v,
\]
where
\[
\beta(r):=|\alpha(r)|^{q(r)-2}\alpha(r).
\]
Therefore
\[
\begin{aligned}
\Delta_{\phi,q}u
&=
\left(\beta e^v\right)_r
+
(n-1)(\log\phi)_r\beta e^v  \\
&=
e^v
\left\{
\beta_r
+
\beta v_r
+
(n-1)(\log\phi)_r\beta
\right\}  \\
&=
e^v
\left\{
\beta_r
+
\beta\left((n-1)(\log\phi)_r+v_r\right)
\right\}.
\end{aligned}
\]
Thus
\[
-\Delta_{\phi,q}u+u^{q(r)-1}
=
e^v
\left\{
-\beta_r
-
\beta\left((n-1)(\log\phi)_r+v_r\right)
+
1
\right\}.
\]
On the other hand, since
\[
f(r,u)=\gamma(r)u^{p(r)-1},
\]
we have
\[
f(r,u)
=
\gamma(r)
\exp\left(
\frac{p(r)-1}{q(r)-1}v
\right).
\]
But
\[
\exp\left(
\frac{p(r)-1}{q(r)-1}v
\right)
=
e^v
\exp\left(
\frac{p(r)-q(r)}{q(r)-1}v
\right).
\]
Therefore
\[
f(r,u)
=
e^v\,
\gamma(r)
\exp\left(
\frac{p(r)-q(r)}{q(r)-1}v
\right)
=
e^v\Gamma(r).
\]
Consequently, the equation
\[
-\Delta_{\phi,q}u+u^{q(r)-1}
=
\gamma(r)u^{p(r)-1}
\]
is equivalent to
\[
e^v
\left\{
-\beta_r
-
\beta\left((n-1)(\log\phi)_r+v_r\right)
+
1
\right\}
=
e^v\Gamma(r).
\]
Since $e^v>0$, this is equivalent to
\[
\beta_r
+
\beta\left((n-1)(\log\phi)_r+v_r\right)
=
1-\Gamma(r),
\]
which is \eqref{q-ODE}.

Finally, multiplying \eqref{q-ODE} by $\phi^{n-1}e^v$, we obtain
\[
\left(
\phi^{n-1}e^v\beta
\right)_r
=
\phi^{n-1}e^v(1-\Gamma(r)).
\]
Using the definition of $\Gamma$, this becomes
\[
\left(
\phi^{n-1}e^v\beta
\right)_r
=
\phi^{n-1}
\left[
e^v
-
\gamma(r)
\exp\left(
\frac{p(r)-1}{q(r)-1}v
\right)
\right],
\]
which is \eqref{eq:q_variable_divergence_form}. The proof is complete.
\end{proof}

\begin{remark}
Equation \eqref{q-ODE} is a nonlinear equation for the logarithmic variable
$v$. In contrast with the constant exponent case, the coefficient
\[
\alpha(r)=\left(\frac{v(r)}{q(r)-1}\right)_r
\]
contains both $v_r$ and $q_r$. Indeed,
\[
\alpha(r)
=
\frac{v_r(r)}{q(r)-1}
-
\frac{v(r)q_r(r)}{(q(r)-1)^2}.
\]
Thus
\[
\beta(r)
=
\left|
\frac{v_r(r)}{q(r)-1}
-
\frac{v(r)q_r(r)}{(q(r)-1)^2}
\right|^{q(r)-2}
\left(
\frac{v_r(r)}{q(r)-1}
-
\frac{v(r)q_r(r)}{(q(r)-1)^2}
\right).
\]
This makes \eqref{q-ODE} substantially more difficult than the constant
exponent equation. For this reason, we regard the variable exponent
reduction mainly as a structural formula rather than as a source of explicit
solutions.
\end{remark}

 \begin{remark}
The reductions in this section provide explicit radial strong solutions and
formal ansatz equations but they are not meant to replace the variational
existence result of Theorem~\ref{reac-diff}, which assures us  the
existence of integrable global $H$-invariant weak solutions. On complete non-compact warped
products, for example on the model
\[
M=\mathcal{R}\times_{\exp(ar)}Y,
\]
we must still verify whether a given explicit radial solution belongs to
$L_1^{q(\cdot)}(M)$. Constant positive solutions, and many of the
non-constant resonant solutions above, fail this integrability condition
because the Riemannian volume is infinite. 

\end{remark} 

\subsection{A decaying manufactured radial solution}

In this subsection we construct explicit positive radial profiles which belong to the
variational space on the $\exp$- and $r$-warped models. Again
\[
\mathcal{R}=[r_0,\infty),\qquad r_0>0,
\]
and
\[
M= \mathcal{R} \times_{\phi}Y,
\qquad
g_\phi=dr^2+\phi(r)^2g_Y,
\]
where $Y$ is a compact $m$-dimensional homogeneous manifold and where
\[
\phi(r)=e^{ar}
\qquad\text{or}\qquad
\phi(r)=ar,
\qquad a>0.
\]
We consider the constant exponent $k>1$ and the equation
\begin{equation}\label{eq:manufactured}
-\Delta_{\phi,k}u+u^{k-1}
=
\gamma(r)u^{k-1}.
\end{equation}
We look for a positive radial solution of the form
\[
u(r)=\phi(r)^{-\sigma},
\qquad
\sigma>0.
\]

\begin{proposition}[A decaying explicit radial profile]
\label{prop:decaying-manufactured-profile}
Let $k\geq2$, $m= n -1 = \dim Y\geq1$, and $a>0$.

\begin{enumerate}
\item If
\[
\phi(r)=e^{ar},
\]
choose
\[
\sigma>\frac{m}{k}
\]
and assume that
\[
1-\sigma^{k-1}a^k\bigl(\sigma(k-1)-m\bigr)>0.
\]
Then
\[
u(r)=e^{-\sigma ar}
\]
belongs to $L^k_1(M)$, and there exists a positive bounded function
$\gamma:\mathcal{R} \to(0,\infty)$ such that $u$ solves
\[
-\Delta_{\phi,k}u+u^{k-1}
=
\gamma(r)u^{k-1}
\]
in the radial strong sense. More explicitly,
\begin{equation}\label{eq:gamma-exp}
\gamma(r)
=
1-\sigma^{k-1}a^k\bigl(\sigma(k-1)-m\bigr).
\end{equation}

\item If
\[
\phi(r)=ar,
\]
choose
\[
\sigma>\frac{m+1}{k}.
\]
Assume, if necessary after increasing $r_0$, that
\[
1-
\sigma^{k-1}\bigl((\sigma+1)(k-1)-m\bigr)r^{-k}
>0
\qquad\text{for every }r\geq r_0.
\]
Then
\[
u(r)=(ar)^{-\sigma}
\]
belongs to $L^k_1(M)$, and there exists a positive bounded function
$\gamma: \mathcal{R}\to(0,\infty)$ such that $u$ solves
\[
-\Delta_{\phi,k}u+u^{k-1}
=
\gamma(r)u^{k-1}
\]
in the radial strong sense. More explicitly,
\begin{equation}\label{eq:gamma-linear}
\gamma(r)
=
1-
\sigma^{k-1}\bigl((\sigma+1)(k-1)-m\bigr)r^{-k}.
\end{equation}
\end{enumerate}
\end{proposition}

\begin{proof}
We treat both cases at once. Set
\[
h(r):=\frac{\phi_r(r)}{\phi(r)}.
\]
Since
\[
u(r)=\phi(r)^{-\sigma},
\]
we have
\[
u_r=-\sigma h u.
\]
Hence
\[
|u_r|^{k-2}u_r
=
-\sigma^{k-1}h^{k-1}u^{k-1}.
\]
Using
\[
u^{k-1}=\phi^{-\sigma(k-1)}
\]
and the radial $k$-Laplacian
\[
\Delta_{\phi,k}u
=
\frac1{\phi^m}
\left(
\phi^m|u_r|^{k-2}u_r
\right)_r,
\]
we obtain
\[
\Delta_{\phi,k}u
=
-\sigma^{k-1}
\frac1{\phi^m}
\left(
\phi^{m-\sigma(k-1)}h^{k-1}
\right)_r.
\]
Therefore
\[
\Delta_{\phi,k}u
=
-\sigma^{k-1}\phi^{-\sigma(k-1)}
\left[
\bigl(m-\sigma(k-1)\bigr)h^k
+
(k-1)h^{k-2}h_r
\right].
\]
Thus
\[
-\Delta_{\phi,k}u+u^{k-1}
=
\phi^{-\sigma(k-1)}
\left[
1+
\sigma^{k-1}
\left(
\bigl(m-\sigma(k-1)\bigr)h^k
+
(k-1)h^{k-2}h_r
\right)
\right].
\]
Since
\[
u^{k-1}=\phi^{-\sigma(k-1)},
\]
we define
\begin{equation}\label{eq:gamma-general}
\gamma(r)
=
1+
\sigma^{k-1}
\left(
\bigl(m-\sigma(k-1)\bigr)h^k
+
(k-1)h^{k-2}h_r
\right).
\end{equation}
Then
\[
-\Delta_{\phi,k}u+u^{k-1}
=
\gamma(r)u^{k-1}.
\]

If $\phi(r)=e^{ar}$, then
\[
h(r)=a,
\qquad
h_r(r)=0.
\]
Hence \eqref{eq:gamma-general} gives
\[
\gamma(r)
=
1+
\sigma^{k-1}
\bigl(m-\sigma(k-1)\bigr)a^k
=
1-\sigma^{k-1}a^k\bigl(\sigma(k-1)-m\bigr),
\]
which is \eqref{eq:gamma-exp}. By assumption this constant is positive.
It is therefore positive and bounded.

We now check that $u\in L^k_1(M)$. Since
\[
dV_{g_\phi}=\phi(r)^m\,dr\wedge dV_{g_Y},
\]
we have
\[
\int_M |u|^k\,dV_{g_\phi}
=
\operatorname{Vol}(Y)
\int_{r_0}^\infty
e^{ar(m-\sigma k)}\,dr.
\]
This integral is finite if and only if
\[
\sigma k>m.
\]
Moreover,
\[
|u_r|^k
=
(\sigma a)^k e^{-\sigma akr},
\]
and therefore
\[
\int_M |u_r|^k\,dV_{g_\phi}
=
(\sigma a)^k\operatorname{Vol}(Y)
\int_{r_0}^\infty
e^{ar(m-\sigma k)}\,dr
<
\infty.
\]
Thus $u\in L^k_1(M)$ in the exponential case.

If $\phi(r)=ar$, then
\[
h(r)=\frac1r,
\qquad
h_r(r)=-\frac1{r^2}.
\]
Therefore \eqref{eq:gamma-general} gives
\[
\gamma(r)
=
1+
\sigma^{k-1}
\left(
m-\sigma(k-1)-(k-1)
\right)r^{-k},
\]
that is,
\[
\gamma(r)
=
1-
\sigma^{k-1}\bigl((\sigma+1)(k-1)-m\bigr)r^{-k},
\]
which is \eqref{eq:gamma-linear}. By the choice of $r_0$, this function
is positive on $\mathcal{R}=[r_0,\infty)$. It is also bounded, since
\[
\gamma(r)\to1
\qquad\text{as }r\to\infty.
\]

It remains to check that $u\in L^k_1(M)$. Since
\[
dV_{g_\phi}=(ar)^m\,dr \wedge\,dV_{g_Y},
\]
we have
\[
\int_M |u|^k\,dV_{g_\phi}
=
\operatorname{Vol}(Y)
\int_{r_0}^\infty
(ar)^{m-\sigma k}\,dr.
\]
This integral is finite if and only if
\[
\sigma k>m+1.
\]
This is exactly ensured by the assumption
\[
\sigma>\frac{m+1}{k}.
\]
Moreover,
\[
u_r
=
-\frac{\sigma}{r}u,
\]
and hence
\[
|u_r|^k
=
\sigma^k r^{-k}(ar)^{-\sigma k}.
\]
Therefore
\[
\int_M |u_r|^k\,dV_{g_\phi}
=
\sigma^k\operatorname{Vol}(Y)
\int_{r_0}^\infty
r^{-k}(ar)^{m-\sigma k}\,dr
<
\infty,
\]
again because $\sigma k>m+1$. Thus $u\in L^k_1(M)$ in the linear case. 
The proof is complete.
\end{proof}

\begin{remark}
We note that
the present
$\exp$- and $r$-warped models and the ansatz
$u=\phi^{-\sigma}$ do not produce a bounded positive
coefficient for superlinear powers $p>k$. Indeed, if we write 
\[
-\Delta_{\phi,k}u+u^{k-1}
=
\gamma(r)u^{p-1},
\qquad p>k,
\]
then the same computation gives
\[
\gamma(r)
=
\phi(r)^{\sigma(p-k)}
\left[
1+
\sigma^{k-1}
\left(
\bigl(m-\sigma(k-1)\bigr)h^k
+
(k-1)h^{k-2}h_r
\right)
\right].
\]
For $\phi(r)=ar$, the bracket tends to $1$, and therefore
$\gamma(r)$ is unbounded as $r\to\infty$. For
$\phi(r)=e^{ar}$, the bracket is constant; if it is positive, then
$\gamma(r)$ is again unbounded for $p>k$, while if it is tuned to
vanish, then $\gamma\equiv0$. Thus the natural bounded positive
manufactured coefficient for these two warpings corresponds to the
case $p=k$.
\end{remark}

\begin{remark}[Compatibility with the variational framework]

The profiles in Proposition~\ref{prop:decaying-manufactured-profile} belong to
$L_1^k(M)$, but the equations with positive bounded manufactured coefficient
have $p=k$. Consequently, they do not satisfy the strict superlinearity
requirements $q\ll p$ and $q^+<p^-$ of Theorem~\ref{reac-diff} when $q=k$.
They are explicit integrable solutions of the resonant equation, rather than
instances of that abstract existence result.
For a superlinear problem with constant exponents to fall within the
subcritical regime of Theorem~\ref{reac-diff}, one would instead need
\[
1<k<m+1,
\qquad
k<p<k^*=\frac{(m+1)k}{m+1-k},
\]
together with all its remaining hypotheses. As shown in the preceding
remark, retaining the same ansatz on the exponential or linear warped end
with $p>k$ does not give a positive bounded manufactured coefficient.
On $[r_0,\infty)\times_\phi Y$, the weak identities for the profiles above
are understood against compactly supported interior test functions; no
boundary condition at $r_0$ is asserted.

\end{remark}

The previous computations provide explicit radial profiles obtained from the one-dimensional reduction of the equation. However, since the original problem is posed on the whole warped product manifold, it is necessary to justify that these radial profiles are not only formal solutions of the reduced equation, but also weak solutions of the full problem.

The following result gives this justification. It shows that, when the radial profile has the required integrability and satisfies the reduced equation in the distributional sense, the corresponding radial function on the full manifold satisfies the weak formulation against arbitrary compactly supported test functions. Therefore, the explicit profiles constructed below can be interpreted as genuine weak solutions of the original equation, provided that they belong to the variational space.

\begin{theorem}[Radial distributional solutions as global weak solutions]
Let $(Y,g_Y)$ be a compact Riemannian manifold of dimension $m\geq 1$, and let
\[
M=\mathcal{R}\times_{\phi}Y,\qquad 
g_{\phi}=dr^2+\phi(r)^2g_Y,
\]
where $\phi\in C^\infty(\mathbb{R})$ and $\phi(r)>0$ for every $r\in\mathbb{R}$. Let $k>1$, $p>1$, and let $\gamma: \mathcal{R} \to\mathbb{R}$ be a measurable function.

Assume that
\[
u\in W^{1,k}_{\mathrm{loc}}(\mathcal{R})
\]
and
\[
\int_{\mathcal{R}}
\left(
|u(r)|^k+|u_r(r)|^k
\right)
\phi(r)^m\,dr<\infty .
\]
Assume also that
\[
\gamma(r)|u(r)|^{p-2}u(r)\in L^1_{\mathrm{loc}}(\mathcal{R})
\]
and that $u$ satisfies
\[
-\frac{1}{\phi(r)^m}
\left(
\phi(r)^m|u_r(r)|^{k-2}u_r(r)
\right)_r
+
|u(r)|^{k-2}u(r)
=
\gamma(r)|u(r)|^{p-2}u(r)
\tag{47}
\]
in the sense of distributions on $\mathcal{R}$. Equivalently, for every
$\eta\in \mathcal{D}(\mathcal{R})$,
\[
\int_{\mathcal{R}}
\phi(r)^m |u_r(r)|^{k-2}u_r(r)\eta_r(r)\,dr
+
\int_{\mathcal{R}}
\phi(r)^m
\left(
|u(r)|^{k-2}u(r)
-
\gamma(r)|u(r)|^{p-2}u(r)
\right)
\eta(r)\,dr
=0 .
\tag{48}
\]

Define $U:M\to\mathbb{R}$ by
\[
U(r,y)=u(r).
\]
Then $U\in L^k_1(M)$, and $U$ is a weak solution on $(M,g_\phi)$ of
\[
-\Delta_{\phi,k}U+|U|^{k-2}U
=
\gamma(r)|U|^{p-2}U .
\tag{49}
\]
That is, for every $\psi\in \mathcal{D}(M)$,
\[
\int_M
\left(
|\nabla_{g_\phi}U|^{k-2}
g_\phi(\nabla_{g_\phi}U,\nabla_{g_\phi}\psi)
+
|U|^{k-2}U\psi
-
\gamma(r)|U|^{p-2}U\psi
\right)
\,dV_{g_\phi}
=0 .
\tag{50}
\]
\end{theorem}

\begin{proof}
Since $U(r,y)=u(r)$, the function $U$ is independent of the variable $y$. In local product coordinates $(r,y^1,\ldots,y^m)$, the inverse metric of
\[
g_\phi=dr^2+\phi(r)^2g_Y
\]
has the form
\[
g_\phi^{-1}
=
\partial_r\otimes\partial_r
+
\phi(r)^{-2}g_Y^{-1}.
\]
Therefore, for every smooth radial function, its Riemannian gradient is given by
\[
\nabla_{g_\phi}U=u_r(r)\partial_r.
\]
The same identity holds for the weak gradient of $U$, because $u\in W^{1,k}_{\mathrm{loc}}(\mathcal{R})$ and the function $U(r,y)=u(r)$ belongs locally to $W^{1,k}$ in product coordinates. Hence
\[
|\nabla_{g_\phi}U|_{g_\phi}=|u_r(r)|
\quad\text{a.e. on }M.
\]

The Riemannian volume form of the warped product is
\[
dV_{g_\phi}
=
\phi(r)^m\,dr \wedge \,dV_{g_Y}.
\]
Since $Y$ is compact, $\operatorname{Vol}_{g_Y}(Y)<\infty$. Thus
\[
\int_M |U|^k\,dV_{g_\phi}
=
\int_{\mathcal{R}}\int_Y
|u(r)|^k\phi(r)^m\,dV_{g_Y}\,dr
=
\operatorname{Vol}_{g_Y}(Y)
\int_{\mathcal{R}}
|u(r)|^k\phi(r)^m\,dr
<\infty.
\]
Similarly,
\[
\int_M |\nabla_{g_\phi}U|^k\,dV_{g_\phi}
=
\int_{\mathcal{R}}\int_Y
|u_r(r)|^k\phi(r)^m\,dV_{g_Y}\,dr
=
\operatorname{Vol}_{g_Y}(Y)
\int_{\mathcal{R}}
|u_r(r)|^k\phi(r)^m\,dr
<\infty.
\]
Consequently, $U\in L^k_1(M)$.

Let now $\psi\in \mathcal{D}(M)$. Define
\[
\Psi(r):=\int_Y \psi(r,y)\,dV_{g_Y}(y).
\]
Since $Y$ is compact and $\psi$ is smooth with compact support in $M$, the function $\Psi$ belongs to $\mathcal{D}(\mathbb{R})$. Moreover, differentiating under the integral sign gives
\[
\Psi_r(r)
=
\int_Y \partial_r\psi(r,y)\,dV_{g_Y}(y).
\]
This is justified because the measure $dV_{g_Y}$ is independent of $r$, and $\psi$ is smooth with compact support.

Set
\[
A(r):=|u_r(r)|^{k-2}u_r(r)
\]
and
\[
B(r):=
|u(r)|^{k-2}u(r)
-
\gamma(r)|u(r)|^{p-2}u(r).
\]
The assumptions imply that $A\in L^1_{\mathrm{loc}}(\mathcal{R})$ and $B\in L^1_{\mathrm{loc}}(\mathcal{R})$. Indeed, on every compact interval $I\subset \mathcal{R}$, the function $\phi$ is bounded above and below by positive constants. Since $u_r\in L^k_{\mathrm{loc}}(\mathcal{R})$, we have
\[
|u_r|^{k-1}\in L^{k/(k-1)}_{\mathrm{loc}}(\mathcal{R})
\subset L^1_{\mathrm{loc}}(\mathcal{R}),
\]
and hence $A\in L^1_{\mathrm{loc}}(\mathcal{R})$. Similarly, $u\in L^k_{\mathrm{loc}}(\mathcal{R})$ implies
\[
|u|^{k-1}\in L^{k/(k-1)}_{\mathrm{loc}}(\mathcal{R})
\subset L^1_{\mathrm{loc}}(\mathcal{R}),
\]
whereas the local integrability of
$\gamma(r)|u(r)|^{p-2}u(r)$ is assumed. Thus all terms in the following identities are well-defined.

By the distributional formulation (48), applied to $\eta=\Psi$, we have
\[
\int_{\mathcal{R}}
\phi(r)^m A(r)\Psi_r(r)\,dr
+
\int_{\mathcal{R}}
\phi(r)^m B(r)\Psi(r)\,dr
=0.
\tag{51}
\]

We now compute the left-hand side of the weak formulation on $M$. Since
\[
\nabla_{g_\phi}U=u_r\partial_r,
\]
we obtain
\[
g_\phi(\nabla_{g_\phi}U,\nabla_{g_\phi}\psi)
=
g_\phi(u_r\partial_r,\nabla_{g_\phi}\psi).
\]
The $r$-component of $\nabla_{g_\phi}\psi$ is $\partial_r\psi$, because the coefficient of $dr^2$ in $g_\phi$ is equal to $1$. Therefore
\[
g_\phi(\nabla_{g_\phi}U,\nabla_{g_\phi}\psi)
=
u_r(r)\partial_r\psi(r,y).
\]
It follows that
\[
|\nabla_{g_\phi}U|^{k-2}
g_\phi(\nabla_{g_\phi}U,\nabla_{g_\phi}\psi)
=
|u_r(r)|^{k-2}u_r(r)\partial_r\psi(r,y)
=
A(r)\partial_r\psi(r,y).
\]

Using the expression of the volume form and Fubini's theorem, we get
\[
\begin{aligned}
&\int_M
\left(
|\nabla_{g_\phi}U|^{k-2}
g_\phi(\nabla_{g_\phi}U,\nabla_{g_\phi}\psi)
+
|U|^{k-2}U\psi
-
\gamma(r)|U|^{p-2}U\psi
\right)
\,dV_{g_\phi}
\\
&=
\int_{\mathcal{R}}\int_Y
\left(
A(r)\partial_r\psi(r,y)
+
B(r)\psi(r,y)
\right)
\phi(r)^m\,dV_{g_Y}(y)\,dr.
\end{aligned}
\]
Since $\psi$ has compact support, and since $A$ and $B$ are locally integrable, the previous integral is finite. We may therefore integrate first in the $Y$-variable. This gives
\[
\begin{aligned}
&\int_M
\left(
|\nabla_{g_\phi}U|^{k-2}
g_\phi(\nabla_{g_\phi}U,\nabla_{g_\phi}\psi)
+
|U|^{k-2}U\psi
-
\gamma(r)|U|^{p-2}U\psi
\right)
\,dV_{g_\phi}
\\
&=
\int_{\mathcal{R}}
\phi(r)^m A(r)
\left(
\int_Y\partial_r\psi(r,y)\,dV_{g_Y}(y)
\right)
dr
\\
&\quad+
\int_{\mathcal{R}}
\phi(r)^m B(r)
\left(
\int_Y\psi(r,y)\,dV_{g_Y}(y)
\right)
dr
\\
&=
\int_{\mathcal{R}}
\phi(r)^m A(r)\Psi_r(r)\,dr
+
\int_{\mathcal{R}}
\phi(r)^m B(r)\Psi(r)\,dr.
\end{aligned}
\]
The last expression is zero by (51). Therefore
\[
\int_M
\left(
|\nabla_{g_\phi}U|^{k-2}
g_\phi(\nabla_{g_\phi}U,\nabla_{g_\phi}\psi)
+
|U|^{k-2}U\psi
-
\gamma(r)|U|^{p-2}U\psi
\right)
\,dV_{g_\phi}
=0
\]
for every $\psi\in \mathcal{D}(M)$. Hence $U$ is a weak solution of (49) on $M$. The proof is complete.
\end{proof}

\begin{corollary}[The decaying radial profile is a global weak solution]
Let $ \mathcal{R} =[r_0,\infty)$, with $r_0>0$, and let
\[
M= \mathcal{R} \times_\phi Y,
\qquad
g_\phi=dr^2+\phi(r)^2g_Y,
\]
where $Y$ is a compact $m$-dimensional homogeneous manifold. Assume that
\[
\phi(r)=e^{ar}
\qquad\text{or}\qquad
\phi(r)=ar,
\qquad a>0.
\]
Let $u=\phi^{-\sigma}$ be one of the profiles constructed in
Proposition~\ref{prop:decaying-manufactured-profile}, and let $\gamma$ be
the corresponding positive bounded coefficient defined by
\eqref{eq:gamma-exp} or \eqref{eq:gamma-linear}. Then
\[
u\in L^k_1(M)
\]
and $u$ is a global weak solution of
\[
-\Delta_{\phi,k}u+u^{k-1}
=
\gamma(r)u^{k-1}.
\]
Equivalently, for every $\psi\in \mathcal{D}(\operatorname{int}M)$,
\[
\int_M
\left(
|\nabla u|^{k-2}g_\phi(\nabla u,\nabla\psi)
+
u^{k-1}\psi
-
\gamma(r)u^{k-1}\psi
\right)
\,dV_{g_\phi}
=0.
\]
\end{corollary}

\begin{proof}
By Proposition~\ref{prop:decaying-manufactured-profile}, the profile
\[
u(r)=\phi(r)^{-\sigma}
\]
belongs to $L^k_1(M)$. Moreover, $u$ satisfies the radial strong equation
\[
-\Delta_{\phi,k}u+u^{k-1}
=
\gamma(r)u^{k-1}
\]
on $I$, where $\gamma$ is positive and bounded.

Let $\psi\in \mathcal{D}(\operatorname{int}M)$. Since $u$ is radial, we have
\[
\nabla u=u_r\partial_r,
\qquad
|\nabla u|=|u_r|.
\]
Therefore
\[
|\nabla u|^{k-2}g_\phi(\nabla u,\nabla\psi)
=
|u_r|^{k-2}u_r\,\psi_r.
\]
Using the warped product measure
\[
dV_{g_\phi}=\phi(r)^m\,dr \wedge\,dV_{g_Y},
\]
we obtain
\[
\int_M
|\nabla u|^{k-2}g_\phi(\nabla u,\nabla\psi)\,dV_{g_\phi}
=
\int_\mathcal{R}\int_Y
|u_r|^{k-2}u_r\,\psi_r\,
\phi(r)^m\,dV_{g_Y}\,dr.
\]
Since $\psi$ has compact support in $\operatorname{int}M$, integration by
parts in the $r$-variable gives no boundary term. Hence
\[
\int_M
|\nabla u|^{k-2}g_\phi(\nabla u,\nabla\psi)\,dV_{g_\phi}
=
-\int_\mathcal{R}\int_Y
\frac{1}{\phi(r)^m}
\left(
\phi(r)^m |u_r|^{k-2}u_r
\right)_r
\psi\,
\phi(r)^m\,dV_{g_Y}\,dr.
\]
That is,
\[
\int_M
|\nabla u|^{k-2}g_\phi(\nabla u,\nabla\psi)\,dV_{g_\phi}
=
-\int_M
\Delta_{\phi,k}u\,\psi\,dV_{g_\phi}.
\]
Using the strong equation, we conclude that
\[
\int_M
\left(
|\nabla u|^{k-2}g_\phi(\nabla u,\nabla\psi)
+
u^{k-1}\psi
-
\gamma(r)u^{k-1}\psi
\right)
\,dV_{g_\phi}
=0.
\]
Thus $u$ is a weak solution in the variational sense.
\end{proof}


\begin{remark}[Dependence on the fibre dimension]
\label{rem:fibre-dimension}
We are aware that, in the final radial constructions, the fibre dimension
is no longer displayed explicitly in the abbreviated profile
$u=\phi^{-\sigma}$ or in the notation $\Delta_{\phi,k}$. This is a feature
of the notation and of the restriction to invariant functions. The
parameter $m=\dim Y=n-1$ remains fixed and is retained throughout the
calculations; the resulting equation is a weighted one-dimensional
reduction of a problem on an $(m+1)$-dimensional manifold.

Indeed, the metric scales each of the $m$ fibre directions by $\phi^2$,
so that its volume density contains the factor $\phi^m$. Although a
radial function $U(r,y)=u(r)$ satisfies
$\nabla_{g_\phi}U=u_r\partial_r$ and
$|\nabla_{g_\phi}U|=|u_r|$, its divergence is computed with this full
volume density. Consequently,
\[
\Delta_{\phi,q}u
=\frac{1}{\phi^m}
\left(\phi^m|u_r|^{q(r)-2}u_r\right)_r
=\left(|u_r|^{q(r)-2}u_r\right)_r
+m\frac{\phi_r}{\phi}|u_r|^{q(r)-2}u_r.
\]
The same factor is present in the logarithmic flux
$\phi^m e^v\beta$ in
\eqref{eq:div_form} and \eqref{eq:q_variable_divergence_form}.
Thus eliminating derivatives in the fibre directions does not eliminate
their contribution to the radial divergence.

The dependence on $m$ is also retained in the manufactured coefficients
\eqref{eq:gamma-exp} and \eqref{eq:gamma-linear}, and in the weighted
Sobolev integrals. In particular,
\[
\int_M\left(|U|^k+|\nabla_{g_\phi}U|^k\right)dV_{g_\phi}
=\operatorname{Vol}_{g_Y}(Y)
\int_{r_0}^{\infty}
\left(|u|^k+|u_r|^k\right)\phi^m\,dr.
\]
For $u=\phi^{-\sigma}$, membership in $L_1^k(M)$ therefore requires
$\sigma k>m$ in the exponential model and $\sigma k>m+1$ in the linear
model, exactly as proved above. The absence of $m$ from the displayed
shape of $u$ does not make these thresholds or the equation independent
of dimension: changing $m$ while keeping $\phi$, $k$, and $\sigma$ fixed
changes $\gamma$ and can change whether $u$ is integrable.

Finally, in the passage to unrestricted weak solutions, fibre integration
replaces an arbitrary test function by
$\Psi(r)=\int_Y\psi(r,y)\,dV_{g_Y}(y)$, while the factor $\phi^m$ remains
in every radial integral. For invariant test functions, only the positive
constant $\operatorname{Vol}_{g_Y}(Y)$ can be cancelled from the weak
identity. Such a cancellation does not remove the dimensional weight.
Likewise, the critical exponent in the abstract existence theorem is
$q^*(r)=(m+1)q(r)/(m+1-q(r))$, under its hypothesis $q^+<m+1$;
it uses the full manifold dimension rather than the dimension of the
radial coordinate.
\end{remark}

\section{Conclusions}\label{sec:conclusions}
We have examined how warped product geometry and compact symmetry enter
nonlinear elliptic equations with variable diffusion exponents. For a
one-dimensional base and a compact homogeneous fibre, invariant
functions depend only on the base coordinate, while the fibre geometry
enters the reduced equation through the volume weight $\phi^m$, with
$m=\dim Y$. The curvature inequalities, volume estimates, and orbit-growth
conditions discussed in Section~\ref{sec:geometry} describe the geometric
requirements relevant to the abstract existence theorem recalled from
\cite{JFA}. When all its geometric and analytical hypotheses hold, that
theorem provides a non-trivial invariant variational solution, and
Corollary~\ref{cor:unrestricted-weak-formulation} gives its weak formulation
against unrestricted test functions.

For positive radial profiles, we established logarithmic reductions in
both the constant and variable exponent settings. In the constant case,
$v=(k-1)\log u$ yields equivalent flux and divergence identities under
the stated regularity assumptions. When $q=q(r)$, the substitution
$v=(q(r)-1)\log u$ exposes the additional dependence on $q_r$ and gives a
structural equation for the transformed profile. These reductions also
clarify the distinction between solving the radial differential equation
and belonging to the weighted Sobolev space: non-zero constant profiles
on the infinite-volume models solve suitable reduced equations but fail
the required integrability condition.

On the expanding ends $[r_0,\infty)\times_\phi Y$, we constructed the
positive decaying profiles $u=e^{-\sigma a r}$ for $\phi=e^{ar}$ and
$u=(ar)^{-\sigma}$ for $\phi=a r$. Their membership in $L_1^k(M)$ is assured by the conditions 
$\sigma k>m$ and $\sigma k>m+1$, respectively, in which $m = n-1$ is the dimension of the fiber $Y$. Under the positivity restrictions of
Proposition~\ref{prop:decaying-manufactured-profile}, these profiles solve equations
with explicitly determined positive bounded coefficients. These are
manufactured solutions in the resonant case $p=k$. For the same ansatz
and $p>k$, the coefficient cannot remain both positive and bounded on
these ends; hence this construction does not realize the strict
superlinear regime of Theorem~\ref{reac-diff}: there is a distinction between the abstract existence theory,
the radial identities, and some explicitly verified resonant solutions.

Finally, integration over the compact fibre establishes that an
integrable radial distributional solution satisfies the weak equation
against arbitrary compactly supported test functions on the full
product. For the half-line models, these tests are supported in the
interior, and no boundary condition at $r_0$ is imposed. Very interestingly, in spite of considering a symmetry-reduced equation, some geometric properties such as 
the fibre dimension and the Ricci scalar curvature, remain present in the operator, the coefficients, and the admissible decay rates. 

\section*{Declarations}

\paragraph{Funding.}
D.J.P and E.G.R. acknowledge support from ANID-Chile via the FONDECYT grant \#1241719. 

\paragraph{Competing interests.}
The authors declare no competing interests.

\paragraph{Data availability.}
This is a theoretical mathematical study. No empirical datasets were generated or analysed. The analytical constructions and derivations
supporting the results are included in the article.

\paragraph{Acknowledgments.} AI tools were used to improve the language and organization of this manuscript;
all mathematical arguments and proofs in the final version were written entirely by the authors.


\end{document}